\documentclass[11pt]{amsart}

\usepackage{amsthm}
\usepackage{amsmath}
\usepackage{amssymb}
\usepackage{graphicx}
\usepackage{enumerate}
\usepackage{mathrsfs}

\usepackage{geometry}
\usepackage{changepage}
\usepackage[numbers]{natbib}
\usepackage{hyperref}
\hypersetup{
pdfpagemode=UseThumbs,
pdftoolbar=true,        
pdfmenubar=true,        
pdffitwindow=false,     
pdfstartview={Fit},    
pdftitle={Surface waves in a channel with thin tunnels and wells at the bottom: non-reflecting underwater topography},    
pdfauthor={},     
pdfsubject={},  
pdfcreator={},   
pdfproducer={}, 
pdfkeywords={}, 
pdfnewwindow=true,      
colorlinks=true,       
linkcolor=magenta,          
citecolor=red,        
filecolor=cyan,      
urlcolor=blue           
}
\allowdisplaybreaks
\numberwithin{equation}{section}
\numberwithin{figure}{section}

\theoremstyle{plain}
\newtheorem{theorem}{Theorem}[section]
\newtheorem{proposition}[theorem]{Proposition}
\newtheorem{lemma}[theorem]{Lemma}
\newtheorem{corollary}[theorem]{Corollary}
\newtheorem{example}[theorem]{Example}
\newtheorem{remark}[theorem]{Remark}
\newtheorem{definition}[theorem]{Definition}
\newtheorem{conjecture}{Conjecture}

\def\BET{\begin{theorem}}
\def\ENT{\end{theorem}}
\def\BEP{\begin{proposition}}
\def\ENP{\end{proposition}}
\def\BEL{\begin{lemma}}
\def\ENL{\end{lemma}}
\def\BEC{\begin{corollary}}
\def\ENC{\end{corollary}}
\def\BEE{\begin{example} \rm}
\def\ENE{\end{example}}
\def\BER{\begin{remark} \rm}
\def\ENR{\end{remark}}
\def\BED{\begin{definition} \rm}
\def\END{\end{definition}}
\def\BECJ{\begin{conjecture}}
\def\ENCJ{\end{conjecture}}

\def\bea{\begin{eqnarray}}
\def\eea{\end{eqnarray}}

\def\beas{\begin{eqnarray*}}
\def\eeas{\end{eqnarray*}}

\def\beq{\begin{equation}}
\def\eeq{\end{equation}}

\def\row{ \nonumber \\ & & }

\def\bff{{\bf f}}

\def\bfu{{\bf u}}

\def\bfA{{\bf A}}

\def\bfR{{\bf R}}

\def\bfW{{\bf W}}

\def\bbB{{\mathbb B}}
\def\bbC{{\mathbb C}}

\def\bbM{{\mathbb M}}
\def\bbN{{\mathbb N}}

\def\bbR{{\mathbb R}}
\def\bbS{{\mathbb S}}

\def\cF{{\mathcal F}}

\def\cH{{\mathcal H}}

\def\cR{{\mathcal R}}

\def\cU{{ U}}
\def\cV{{\mathcal V}}
\def\cW{{\mathcal W}}
\def\cX{{\mathcal X}}
\def\cY{{\mathcal Y}}

\def\ef{\eqref}

\def\sfk{{\sf k}}
\newcommand{\eps}{\varepsilon}
\newcommand{\dsp}{\displaystyle}
\newcommand{\mrm}[1]{\mathrm{#1}}

\begin{document}

\title[Invisible topography]{Surface waves in a channel with thin tunnels and wells\\[2pt]at the bottom: non-reflecting underwater topography}

\author[L. Chesnel]{Lucas Chesnel}
\author[S. A. Nazarov]{Sergei A. Nazarov}
\author[J. Taskinen]{Jari Taskinen}

\address{INRIA/Centre de math\'ematiques appliqu\'ees, \'Ecole Polytechnique, Universit\'e Paris-Saclay, Route de Saclay, 91128 Palaiseau, France.}
\email{lucas.chesnel@inria.fr}

\address{Saint-Petersburg State University,
Universitetskaya nab., 7--9,  St. Petersburg, 199034, Russia, and  
\\ 
Institute of Problems of Mechanical Engineering RAS,
V.O., Bolshoi pr., 61, St. Petersburg, 199178, Russia.}
\email{srgnazarov@yahoo.co.uk}

\address{ Department of Mathematics and Statistics, University of Helsinki,
P.O.Box 68, FI-00014 Helsinki, Finland.} 
\email{jari.taskinen@helsinki.fi}

\maketitle

\begin{center}
\begin{minipage}{0.9\textwidth}
\noindent\textbf{Abstract.} We consider the propagation of surface water waves in a straight planar channel perturbed at the bottom by several thin curved tunnels and wells. We propose a method to construct non reflecting underwater topographies of this type at an arbitrary prescribed wave number. To proceed, we compute asymptotic expansions of the diffraction solutions with respect to the small parameter of the geometry taking into account the existence of boundary layer phenomena. We establish error estimates to validate the expansions using advances techniques of weighted spaces with detached asymptotics. In the process, we show the absence of trapped surface waves for perturbations small enough. This analysis furnishes asymptotic formulas for the scattering matrix and we use them to determine underwater topographies which are non-reflecting. Theoretical and numerical examples are given. \\

\end{minipage}
\end{center}



\section{Introduction.}
\label{sec1}

\subsection{Non-reflecting and invisible obstacles in waveguides.}\label{sec1.3}

\begin{figure}[!ht]
\centering
\includegraphics[width=10cm,trim={0 0.3cm 0 0.4cm},clip]{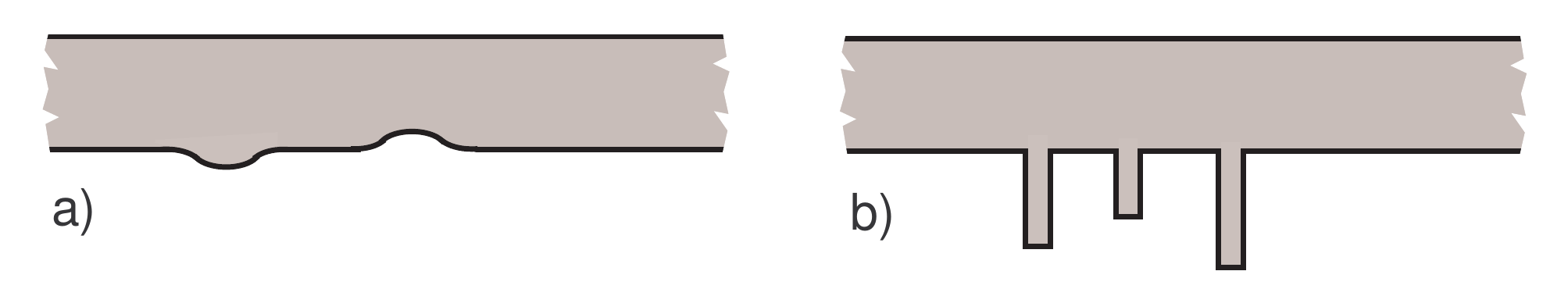}
\caption{a) Gently sloping bottom. b) Bottom with wells.\label{fig5}}
\end{figure}

We investigate the propagation of surface water-waves in a planar channel in time-harmonic regime. We assume that the channel coincides, outside a region where the bottom is geometrically perturbed, with the reference straight channel. We consider a situation where an incident wave propagates through the channel, hits the geometrical defect and gives birth to a scattered field. One commonly denotes by $R$ the \textit{reflection} coefficient, which corresponds to the amplitude of the backscattered farfield, and by $T$ the \textit{transmission} coefficient, which corresponds to the amplitude of the transmitted farfield. Due to conservation of energy, these two complex numbers satisfy the relation
\begin{equation}\label{ConservationNRJ}
|R|^2+|T|^2=1.
\end{equation}
The scattering coefficients $R$, $T$ depend on the geometry and satisfy $R=0$, $T=1$ in the reference straight channel. In this context, a question of growing interest is to find situations where one has good transmission properties (see in particular the literature concerning so-called Perfect Transmission Resonances (PTRs) \cite{Shao94,PoGP99,LeKi01,Zhuk10,MrMK11}). In particular, one can wish to have perfect transmission in energy, that is $|T|=1$. Due to \ef{ConservationNRJ}, this is equivalent to have $R=0$. In this case, following the terminology introduced in \cite{BN}, we shall say that the perturbation of the bottom is \textit{non reflecting}. Note that in this situation, the transmitted wave in general exhibits a phase shift with respect to the incident field. One can be more demanding and look for channels where $R=0$ and $T=1$. In this case, we shall speak of \textit{perfect invisibility}.\\
\newline
Usually $R$ and $T$ depend analytically on the geometric parameters defining the channel/waveguide. As a consequence, non reflecting and invisibility situations are unstable: a small change of the setting may ruin them. Eigenvalues embedded into the continuous spectrum behave in a similar way, see for example  \cite{AsPV00}. In \cite{na489,na546}, the notion of {\it enforced stability} of eigenvalues was introduced and the method of {\it fine-tuning the geometric parameters} to maintain the eigenvalues in the continuous spectrum while perturbing the problem was developed and rigorously substantiated. To summarize, this method boils down to mimic the proof of the implicit functions theorem considering a certain indicator of existence of eigenvalues as a function of the parameters of the problem (geometry of the waveguide, physical coefficient in the equation, ...). It was adapted in \cite{BN} to the problem of invisibility. More precisely, in \cite{BN} the authors study an acoustic problem in a waveguide with locally gently sloping walls (Fig.\,\ref{fig5},\,a)) and propose a method to construct non reflecting perturbations of the reference straight geometry. We emphasize that the approach of \cite{BN} does not allow  one to control the phase shift between the incident field and the transmitted field to impose $T=1$. As a consequence, in general the distortion of the wall is not perfectly invisible.\\
\newline
The method of \cite{BN} was then adapted in \cite{na582} to the water-wave
problem in a planar channel with a gently sloping bottom (see again Fig.\,\ref{fig5},\,a)). In \cite{na582}, it is shown that for the water-wave problem, the above technique does not only allow one to find non straight channels where $R=0$, $|T|=1$ but also geometries where $R=0$, $T=1$ (without phase shift between the incident and the transmitted fields). In other words, it allows one to construct perfectly invisible perturbations of the bottom. To be exhaustive, we must mention however that the technique does not apply for the particular case where, with the notation below, $a^2 \not= 2 \sfk^2$ or equivalently, $\lambda \not= \sqrt{2} \sfk \tanh(\sqrt{2} \sfk d)$.

\noindent In order to obtain perfect invisibility ($R=0$, $T=1$) for the acoustic problem considered in \cite{BN}, another way to perturb the reference straight geometry was studied in \cite{na648}. It consists in working with waveguides with several thin rectangles perpendicular to the wall, see Fig.\,\ref{fig5},\,b). This singular perturbation gives a different differential of the scattering coefficients with respect to the geometry (have in mind the proof of the implicit functions theorem mentioned above) and in \cite{na648}, it is proved that it allows one to achieve perfect invisibility in monomode regime by choosing the dispositions and lengths of the rectangles properly.\\
\newline
In the present work, we adapt the approach of  \cite{na648} to the water-wave problem \ef{3}--\ef{5} in the much more tangled geometry with curved thin tunnels. We will show that the approach of \cite{na648} can be used to get $R=0$ for all $\sfk \geq 0$ and $\lambda > \lambda_\dagger(\sfk)$ (again see the notation below) without the above restriction $a^2 \not= 2 \sfk^2$ which appears when considering smooth perturbations (Fig.\,\ref{fig5},\,a)) of the geometry. We will also prove that this manner of deforming the bottom  does not allow one to get $T=1$. In most of the article, we will perturb the reference straight channel by digging thin tunnels. We emphasize that the results of this article cover the more simple case of rectangular well shaped perturbations, see Fig.\,\ref{fig5},\,b) and \S\ref{sec5.4}.

\subsection{Junctions of domains with different limit dimensions.}
\label{sec1.4} 
\begin{figure}[!ht]
\centering
\includegraphics[width=10cm,trim={0 0.2cm 0 0.5cm},clip]{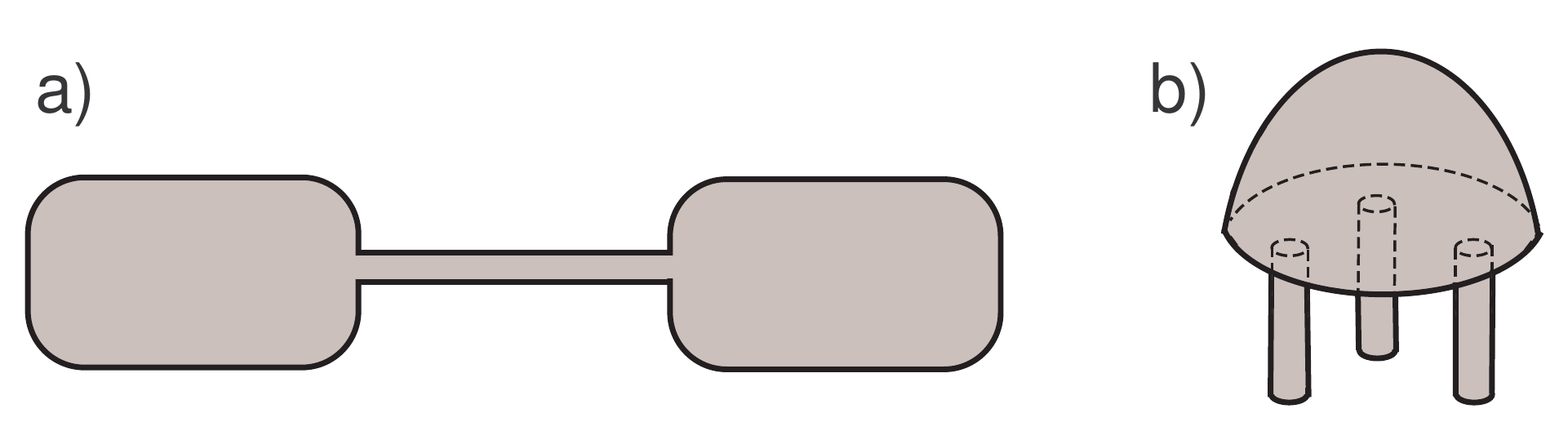}
\caption{a) Dumbbell. b) Junction of a massive body and thin rods.\label{fig6}}
\end{figure}

Invisibility questions are one motivation of the present article. The technique that we will propose to get $R=0$ relies heavily on asymptotic analysis in junctions of massive domains with thin ligaments and developing new results in this field is another goal of the paper. The {\it dumbbell} (see Fig.\,\ref{fig6},\,a)), namely two massive domains $\Omega_\pm$ connected by a thin cylinder $Q^\eps$ with a cross-section of diameter $\eps$, is a classical object in asymptotic analysis. There are many studies of the Neumann Laplacian showing that its spectrum 
\bea
0 = \lambda_1^\eps < \lambda_2^\eps \leq 
\lambda_3^\eps \leq \ldots \leq  \lambda_n^\eps \leq \ldots
\to + \infty \label{ph6}
\eea
has the following distinguishing feature: the limit set 
\bea\label{ph7}
\{\,\lambda_n^0 = \lim\limits_{\eps \to 0^+} \lambda_n^\eps \,\} 
\eea
is the union of three spectra, namely the spectra of the Neumann Laplacian
in $\Omega_\pm$ and the Dirichlet problem for the operator $- \partial_z^2$ 
on the axis of $Q^\eps$. After the pioneering works \cite{Beal73,Arse76}, this
problem has been investigated  in many papers in the original, Fig.\,\ref{fig6},\,a), and modified, Fig.\,\ref{fig6},\,b), formulations with different methods and goals, see e.g. \cite{Gady93, KoMM94,Naza96, na285, na328, A24,Naza05, A30,A31, A32, A4}. 
In particular, complete asymptotic expansions of eigenvalues and eigenfunctions were constructed in \cite{A24} by the methods of  matched asymptotic expansions. Expansions for the  solutions to stationary mixed boundary-value problems were found in \cite{KoMM94,Naza96,na285} by the method of compound asymptotic expansions in dimensions $d \geq3 $ and $d =2$. \\

\noindent It is remarkable that convergence theorems for the spectra \ef{ph6} and \ef{ph7}
have been obtained for domains $\Omega_\pm$ with arbitrary shapes (see e.g.
\cite{Beal73}), but the applications of asymptotic analysis have only been made under additional simplifying assumptions that the boundaries $\partial \Omega_\pm$ are flat near the junctions and that the ligaments $Q^\eps$ are perpendicular to them. The main asymptotic terms are the same for curved surfaces, but the higher-order terms are influenced
by the curvature (cf. \cite{Naza96,na285}). Besides, the derivation of asymptotically sharp estimates becomes much more complicated in the curved case, see \cite{na328,A4} for details. However, the ligaments connecting massive domains are always assumed to be straight.\\
\newline
In this work, we diverge from the traditional formulation for junction problems by, first, considering curved tunnels \ef{1} with variable widths (treating the straight wells Fig.\,\ref{fig5},\,b) as a special case in Section \ref{sec5.4}) and, second, allowing the mid-line of the tunnels to meet the bottom non perpendicularly. On the other hand, to simplify the treatment of the boundary layers, we make the assumption that the boundary of the geometry consists of straight segments near the junction zones.\\
\newline
The first asymptotic terms, which are needed in the fine-tuning procedure to find geometries where $R=0$, are constructed by means of the method of matched asymptotic expansions. However, we emphasize that if one wishes to construct infinite asymptotic series, it is better to apply the method of compound expansions which crucially simplifies the iterative process (cf. \cite{KoMM94,Naza96,na285} and others). The reason is that the limit problems are solved in the same function spaces, whereas the method of matched asymptotic expansions requires for the solutions with singularities of ever growing orders. In order to shorten the article, we decide not to present the construction of infinite asymptotic series. But this can be done.\\
\newline
The very novelty  of the asymptotic analysis in this paper lies in the justification scheme in Section \ref{SectionJustif}. First, we will work with the traditional weighted spaces with detached asymptotics, the norms of which contain the moduli of the scattering coefficients. Consequently, proving error estimates with these norms directly implies the justification of the asymptotics of the scattering matrix. Second, we will consider a Sobolev space endowed with a rather exotic norm, which is defined as the infimum of the norms of several components related to the structure of the junctions (see Section \ref{SectionJustif}). On this occasion, it is worth to mention that the Sobolev and H\"older norms, even their weighted variants with diversified weights on functions themselves and their derivatives in different directions, cannot properly reflect entangled composite asymptotic structures of solutions in junctions of thin and massive domains, especially in the vectorial case, like in elasticity. In this way, our innovative trick allows us to take into account miscellaneous contributions of all geometric parts  of the junctions in the norm. This approach is entirely new and is certainly expected to be helpful in examining other boundary value problems with singular perturbations for example in hybrid domains \cite{Lion06,Naza96,na345,na576} and for elastic junctions \cite{KMM2,KoMR01,Naza05,na514} where reactions of thin fragments on longitudinal and transversal loadings are very discrepant.

\subsection{Outline of the paper.} \label{sec1.5}

We start in Section \ref{SectionSetting} by presenting the setting of the problem. In Section \ref{SectionMainAnalysis}, first we give the main terms  appearing in the asymptotic expansions of the scattering coefficients $R$, $T$. Then we present the fine-tuning procedure introduced in \cite{BN,na582,na648} which allows us to construct non-reflecting underwater topographies for surface waves at any prescribed $\sfk \geq 0$ and $\lambda > \lambda_\dagger(\sfk)$. In Section \ref{SectionNumerics}, we explain how to implement the method numerically and give several examples of non reflecting channels. Section \ref{SectionAsymptotics} is dedicated to the formal asymptotic expansion of the scattering solutions. This provides us in a rather direct way the main asymptotic terms in the decomposition of the scattering coefficients $R^{\eps}_-$, $T^{\eps}$ (with the notation below). The most technical part, Section \ref{SectionJustif}, contains the operator formulation of the problem \ef{3}--\ef{5} and the derivation of asymptotically sharp estimates for the solutions with respect to the norms of weighted spaces with detached asymptotics. We emphasize that these norms are closely connected with the asymptotic structures derived in Section \ref{SectionAsymptotics}, which makes it quite easy to obtain error estimates. In this part, we also prove the absence of trapped modes in the channel $\Pi^\eps$ for $\eps > 0$ small enough. We end the article with some concluding remarks. The main results of this article are Theorem \ref{thNoRe} (non-reflecting geometries) and the approach of Section \ref{SectionJustif} to prove error estimates for asymptotic expansions using well-chosen norms.

\section{Setting of the water-wave problem.}\label{SectionSetting} 

\subsection{Notation.} 
\label{sec1.1}
\begin{figure}[!ht]
\centering
\includegraphics[width=15cm]{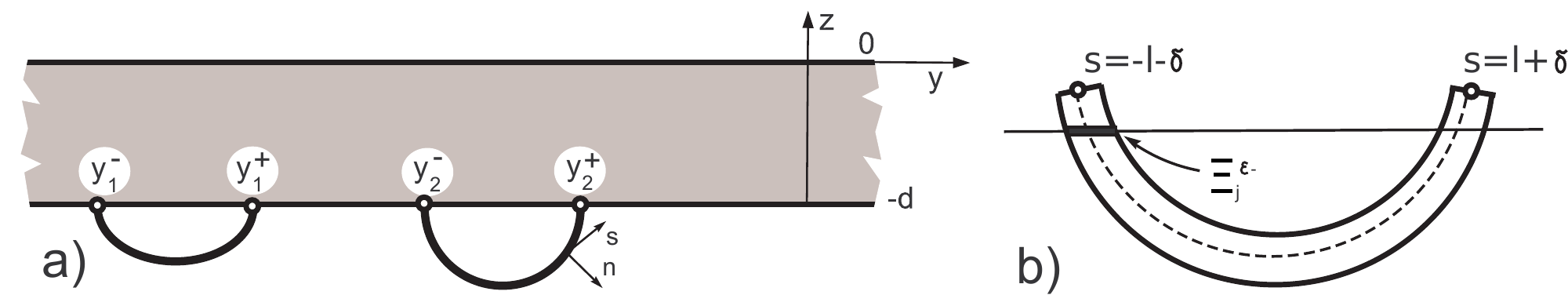}
\caption{a) Straight channel $\Pi$ and defining curves $L_j$. b) Tunnel defined by a curve.\label{fig1}}
\end{figure}

Let 
$\Pi := \{ x= (y,z)\in\bbR\times(-d;0)\}$ be a straight two-dimensional channel. Let $L_j$, $j=1, \ldots, J$ with  $J \in\bbN:=\{0,1,\dots\}$, be a simple smooth curve inside the lower half-plane $\{ x : z < -d\}$, 
connecting the points
\bea
P_j^- = (y_j^-,-d) \qquad\mbox{and} \qquad P_j^+  = (y_j^+,-d)
\qquad\mbox{with} \qquad y_{j-1}^+ < y_j^- < y_j^+ < y_{j+1}^- ,
\label{0}
\eea
see Fig.\,\ref{fig1},\,a). The length of $L_j$ is $2\ell_j > 0$. In a neighbourhood of $L_j$ we introduce the local coordinates $(n_j,s_j)$ where $s_j \in (-\ell_j;\ell_j)$ is the arc length and 
$n_j$ is the oriented distance to $L_j$. We assume that $L_j$
intersects the bottom $\Gamma_d := \bbR\times\{-d\}$ of
the channel at the angles $\alpha_j^\pm \in (0;\pi)$ (with the line $(y_j^{\pm};+\infty)\times\{-d\}$) and denote by 
$\widehat L_j$ a  smooth extension of $L_j$ inside $\Pi$ for the values
$s \in [ -\ell_j -\delta_j;\ell_j + \delta_j]$ for some
small $\delta_j > 0$. We define the domain 
\bea
\mathcal{T}_j^\eps := \big\{ x : s_j \in ( -\ell_j -\delta_j;\ell_j + 
\delta_j), \ -\eps H_j^- (s_j) < n_j < \eps H_j^+ (s_j) \big\}
\label{1}
\eea
entering the channel $\Pi$ (see Fig.\,\ref{fig1},\,b)) and set the thin curved strip (tunnel) $\varpi_j^\eps :=\mathcal{T}_j^\eps\setminus\Pi$. Here $H_j^\pm \in\mathscr{C}^\infty[- \ell_j -\delta_j;\ell_j + \delta_j]$ are  smooth profile functions such that $H_j := H_j^+ + H_j^- > 0$ and $\eps > 0$ is a small parameter.

\begin{figure}[!ht]
\centering
\includegraphics[width=8cm]{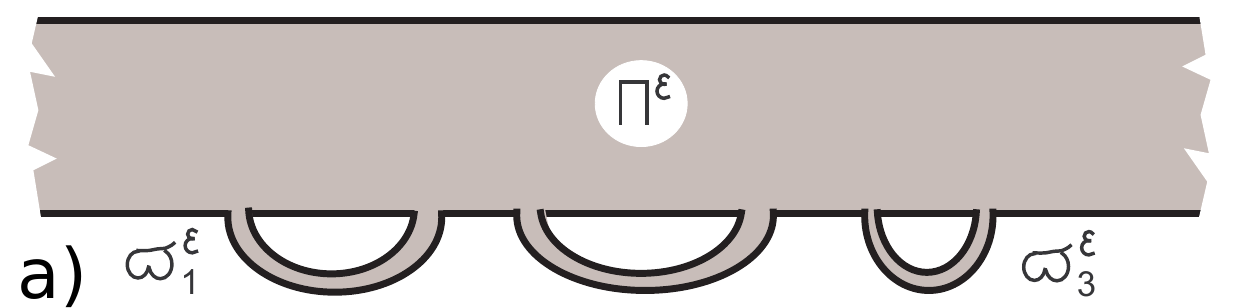}
\includegraphics[width=8cm]{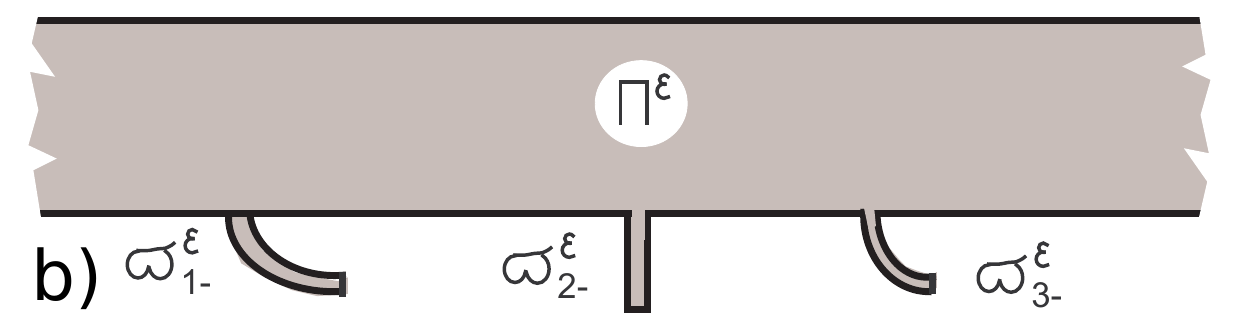}
\caption{a) Tunnels and b) wells in the water domain $\Pi^\eps$.\label{fig2}}
\end{figure}

\noindent We consider the channel
\bea
\Pi^\eps := \Pi \cup \varpi_1^\eps \cup \ldots \cup \varpi_J^\eps \label{2}
\eea
with several thin curved tunnels under the bottom, see Fig.\,\ref{fig2},\, a). We denote $\Gamma_0 := \bbR\times\{0\}$ the free surface of $\Pi^{\eps}$. We study the water-wave problem consisting of the Helmholtz equation
\bea
- \Delta u^\eps + {\sf k}^2 u^\eps =  0 \ \mbox{ in }\Pi^\eps,
\label{3}
\eea
the Neumann (no penetration) boundary condition
\bea
\partial_\nu u^\eps = 0\ \mbox{ on }\partial \Pi^\eps
\setminus \Gamma_0, \label{4}
\eea
and the (kinematic) Steklov condition on the free surface 
\bea
\partial_\nu u^\eps = \lambda u^\eps\ \mbox{ on }\Gamma_0 . \label{5}
\eea
Here, $u^\eps$ is the velocity potential, $\Delta$ is the Laplace
operator and $\partial_\nu$ stands for the outward normal derivative so that $\partial_\nu=
\partial_z$ on $\Gamma_0$. Moreover ${\sf k} \geq0$ is the wave number in the direction perpendicular to the plane $\bbR^2 \ni x$ while $\lambda = g^{-1} \omega^2$ is the spectral parameter, $\omega>0$ being the frequency of time-harmonic oscillations and $g>0$ the acceleration of gravity. Note that in \S\ref{sec5.4}, we will consider wells as drawn in Fig.\,\ref{fig2} ,\,b). In this case, we impose the Neumann condition (\ref{4}) at the ends of the wells.

\subsection{Surface waves and scattering matrix.} 
\label{sec1.2} 
It is known, see e.g. \cite{Have29,KuMaVa}, that the 
continuous spectrum of the problem \eqref{3}--\eqref{5} coincides with the closed
semi-axis $[\lambda_\dagger;+\infty)$, where the cut-off value is given
by
\[
\lambda_\dagger = \lambda_\dagger (\sfk) := \sfk \tanh (\sfk d)=\sfk \frac{e^{\sfk d} -e^{-\sfk d}}{
e^{\sfk d} + e^{-\sfk d}}  . 
\]
Note in particular that for $\sfk=0$, we have $\lambda_\dagger=0$. For any $\lambda > \lambda_\dagger$, introduce the functions $w^\pm$ such that
\bea
w^\pm (y,z) = e^{\pm i y\sqrt{a^2 - \sfk^2} }
\big( e^{a(z+d) } + e^{- a(z+d) } \big) , \label{7}
\eea
where the exponent $a = a( \lambda)>0$ solves the equation 
\[
a \tanh (ad)=\lambda. 
\]
Observe that $w^\pm$ satisfy the Helmholtz equation (\ref{3}), 
the Neumann condition (\ref{4}) on the bottom $\Gamma_d$ of the straight channel and the Steklov condition on the free surface $\Gamma_0$. They correspond to surface waves propagating along the channel $\Pi$ from $\mp \infty$ to $\pm \infty$ (with a convention of a time-harmonic regime in $e^{-i\omega t}$).\\
\newline
Now we focus our attention on the scattering of the waves $w^\pm$ in the singularly perturbed channel $\Pi^\eps$. To describe separately the behaviours at $+\infty$ and $-\infty$, first we introduce two smooth cut-off functions $\chi_\pm \in \mathscr{C}^\infty(\bbR)$ such that $0 \leq \chi_\pm \leq 1$,
\begin{equation}\label{CutOffChipm}
\chi_\pm  = 1 \quad\mbox{for} \ \pm y \geq 2 \ell_0 \quad\qquad\mbox{and}\quad\qquad
\chi_\pm = 0 \quad\mbox{for} \ \pm y \leq  \ell_0 .
\end{equation}
The parameter $\ell_0>0$ is chosen such that the tunnels \ef{1} are all contained in the region $\{ x :|y| < \ell_0 \}$ (in particular, we have $- \ell_0 < y_1^- < y_1^+ < \ldots < y_J^- < y_J^+ < \ell_0$). In the following, we shall say that a function $u$ satisfying \ef{3}--\ef{5} is \textit{outgoing} if it admits the expansion
\begin{equation}\label{DefOut}
u=\chi_-\,A\,w^-+\chi_+\,B\,w^++\tilde{u}
\end{equation}
for some coefficients $A$, $B\in\bbC$ and some remainder $\tilde{u}$ which decays exponentially for $|y|\ge \ell_0$. More precisely, one can prove that necessarily $\tilde{u}$ decays with the rate $O( e^{- \sqrt{b^2 + \sfk^2 }|y|}) $, where $b$ is the root of the transcendental equation 
\begin{equation}\label{TransEqb}
-b \tan (bd)=\lambda
\end{equation}
in $(\pi/(2d),3\pi/(2d))$. One can show that the problem \ef{3}--\ef{5} in $\Pi^\eps$ admits a solution $u_{-}^\eps$ (resp. $u_{+}^\eps$) such that $u_{-}^\eps-w^{+}$ (resp. $u_{+}^\eps-w^{-}$) is outgoing. Due to \ef{DefOut},  we have the representations \begin{equation}\label{totalFields}
\begin{array}{l}
u_-^\eps=\chi_-\,(w^++R^{\eps}_-\,w^-)+\chi_+\,T^{\eps}_-w^++\tilde{u}_-^\eps\\[5pt]
u_+^\eps=\chi_+\,(w^-+R^{\eps}_+\,w^+)+\chi_-\,T^{\eps}_+\,w^-+\tilde{u}_+^\eps,
\end{array}
\end{equation}
where $R^{\eps}_{\mp}$, $T^{\eps}_{\mp}\in\bbC$ and $\tilde{u}_{\mp}^\eps$ decay exponentially for $|x|\ge \ell_0$. The function $u_{\mp}^\eps$ represents the total field associated with the incident field $w^{\pm}$ incoming from $\mp\infty$. Uniqueness for $u_{\mp}^\eps$ occurs if and only if trapped modes are absent. We remind the reader that trapped modes are solutions to \ef{3}--\ef{5} (without source term nor incident field) which decay exponentially for $|x|\ge \ell_0$. They can appear only for  a discrete set of wavenumbers $\sfk$ (see e.g. \cite{EvLV94}). In what follows (see Theorem \ref{th4E}), we shall prove that trapped modes do not exist for $\eps$ small enough so that $u_{\mp}^\eps$ are well-defined.  In (\ref{totalFields}), $R^{\eps}_{\mp}$ and $T^{\eps}_{\mp}$ are usually called \textit{reflection} and \textit{transmission} coefficients. It is known that $T^{\eps}_{-}=T^{\eps}_{+}$ and to simplify we shall denote $T^{\eps}=T^{\eps}_{-}=T^{\eps}_{+}$. The scattering matrix 
\begin{equation}\label{S1}
\bbS^\eps := 
\left( 
\begin{array}{cc} 
R_{-}^\eps & T^\eps \\
T^\eps & R_{+}^\eps 
\end{array}
\right)
\end{equation} 
is uniquely defined, symmetric and unitary. In particular, we have 
\begin{equation}\label{conservationNRJ}
|R_{-}^\eps|^2 + |T^\eps|^2
= 1 \qquad\mbox{ and }\qquad |R_{+}^\eps|^2 + |T^\eps|^2
= 1\qquad\mbox{(conservation of energy).}
\end{equation}
When $R_{-}^\eps=0$, the backscattered field associated with $u^{\eps}_-$ (see (\ref{totalFields})) is evanescent. Note that in this situation, due to (\ref{conservationNRJ}), we also have $R_{+}^\eps=0$, $|T^{\eps}|=1$ and as mentioned in the introduction, we say that the family of tunnels $\varpi_1^\eps, \ldots, \varpi_J^\eps$ is non-reflecting.

\section{Non-reflecting topographies}\label{SectionMainAnalysis}
 Below we will compute asymptotic expansions of the functions $u_{\pm}^{\eps}$ defined in (\ref{totalFields}) with respect to $\eps$. Then we will derive expansions of the scattering coefficients of the matrix $\mathbb{S}^{\eps}$. Since the procedure is a bit long, we first give the main results and explain how to use them to construct non-reflecting topographies, that is to obtain $R^{\eps}_-=0$. For the construction of the asymptotics and the proof of error estimates, we refer the reader to Sections \ref{SectionAsymptotics} and \ref{SectionJustif} respectively.

\subsection{Asymptotics of the scattering coefficients.}  In what follows, for the reflection and transmission coefficients we shall consider the expansions:
\begin{equation}\label{MainAnsatzRef}
R^{\eps}_-=R^{0}_-+\eps R^{'}_-+\eps\widetilde{R}^{\eps}_-;\qquad\qquad T^{\eps}=T^{0}+\eps T^{'}+\eps\widetilde{T}^{\eps}.
\end{equation}
Here $R^{0}_-$, $R^{'}_-$, $T^{0}$, $T^{'}$ are complex constants which are independent of $\eps$ and  $\widetilde{R}^{\eps}_-$, $\widetilde{T}^{\eps}$ correspond to some abstract remainders. First we will show that $R^{0}_-=0$ and $T^{0}=1$. This is natural because when $\eps\to0$, the thin tunnels disappear and the incident wave $w^+$ propagates in $\Pi$ without being perturbed. Next for the correction terms in (\ref{MainAnsatzRef}), we will establish the following important formulas
\begin{equation}\label{MainCorR}
\phantom{\mbox{and }}\hspace{1.2cm} R'_-=\cfrac{-i}{2N\sqrt{a^2-k^2}}\,\dsp\sum_{j=1}^J\int_{-\ell_j}^{\ell_j} H_j(s)\left(\bigg(\cfrac{dv^0_j}{ds}(s)\bigg)^2+\sfk^2(v^0_j(s))^2\right)\,ds
\end{equation}
\begin{equation}\label{MainCorT}
\mbox{and }\qquad T'=\cfrac{-i}{2N\sqrt{a^2-k^2}}\,\dsp\sum_{j=1}^J\int_{-\ell_j}^{\ell_j} H_j(s)\left(\bigg|\cfrac{dv^0_j}{ds}(s)\bigg|^2+\sfk^2|v^0_j(s)|^2\right)\,ds.
\end{equation}
Here $N:=(2a)^{-1}\,(e^{2ad}-e^{-2ad})+2d$ is a normalisation factor and $H_j(s) = H^+_j(s) + H^-_j(s) > 0$ corresponds to the rescaled thickness of the curved strip $\varpi^\eps_j$. Moreover, the functions $v^0_j$  are defined as the solutions to the problems
\begin{equation}\label{Problem1D}
\begin{array}{|lcl}
\ -\cfrac{d}{ds}\Big(H_j(s)\cfrac{dv^0_j}{ds} (s)\Big)+\sfk^2H_j(s)v^0_j(s)=0\qquad\mbox{for }s\in(-\ell_j;\ell_j)\\[12pt]
\hspace{3cm}v^0_j(\pm \ell_j)=w^{+}(y^\pm_j,-d).
\end{array}
\end{equation}
In \S\ref{paragraphCalculus}, we will give explicit examples of setting where the geometrical parameters are such that $R'_-=0$. In this case, the family of tunnels is almost non-reflecting. More precisely, in this situation a perturbation of the straight channel $\Pi$ of order $\eps$ produces a reflection of order $\eps^2$ only. However, due to the presence of the remainder $\widetilde{R}^{\eps}_-$ in the representation (\ref{MainAnsatzRef}), the identity $R'_-=0$ does not yet suffice to guarantee the non-reflectability. Therefore, we have to refine our strategy.
\subsection{The fine-tuning procedure.\label{FineTunningParag}} Let us fix once for all a family of geometrical parameters 
\bea\label{R2}
y_j^\pm, L_j, \ell_j, H_j^\pm, \ \ j= 1, \ldots, J, 
\eea
such that $R'_-=0$. Now we shift two tunnels along the abscissa axis. More precisely, for $h=(h_1,h_2)^{\top}\in \bbR^2$ we assume that $\varpi_1^\eps$, $\varpi_2^\eps$ (reindex the tunnels if necessary) are changed into $\varpi_1^{\eps}(h)$, $\varpi_2^{\eps}(h)$ with 
\bea
\varpi_j^{\eps}(h):= \{ x= (y,z) : (y -h_j, z) \in \varpi_j^\eps \}, \ \ j=1,2.
\label{R3}
\eea
We regard $h=(h_1,h_2)^{\top}\in \bbR^2$ as small parameters which are independent of $\eps$. We denote respectively $R^{\eps}_-(h)$ and $T^{\eps}(h)$ the reflection and transmission coefficients in the geometry 
\begin{equation}\label{ModifiedChannel}
\Pi^\eps(h) := \Pi \cup \varpi_1^\eps(h) \cup\varpi_2^\eps(h) \cup \varpi_3^\eps\cup \ldots \cup \varpi_J^\eps .
\end{equation}
Similarly to (\ref{MainAnsatzRef}), we have the expansion 
\begin{equation}\label{ExpansionEpsh}
R^{\eps}_-(h)=0+\eps R_-'(h)+\eps\widetilde{R}^{\eps}_-(h).
\end{equation}
Here $R'_-(h)$ is given by (\ref{MainCorR}) with, for $j=1,2$, $v^0_j$ defined as the solution to (\ref{Problem1D}) with a data $w^{+}(y^\pm_j,-d)$ replaced by $w^{+}(y^\pm_j+h_j,-d)$. One observes that the map $h\mapsto R'_-(h)$ is analytic. Therefore, we have the expansion
\begin{equation}\label{Expansionh}
R_-'(h)=R_-'(0,0)+h_1\cfrac{\partial R_-'}{\partial h_1}(0,0)+h_2\cfrac{\partial R_-'}{\partial h_2}(0,0)+\widetilde{R}_-'(h)
\end{equation}
where $\widetilde{R}_-'(h)$ is an abstract remainder. Using that $R_-'(0,0)=0$ (this results from the particular choice of the geometrical parameters (\ref{R2})), we find that there holds $R^{\eps}_-(h)=0$ if and only $h\in\bbR^2$ solves the following system of transcendental equations
\begin{equation}
\begin{array}{|lcl}
\  0 &= &\dsp\eps^{-1} \Re e\,R_-^\eps(h) = h_1\,\Re e\, \frac{\partial R_-'}{\partial h_1}	(0,0)+ h_2\,\Re e\, \frac{\partial R_-'}{\partial h_2 }	(0,0) + \Re e \,\widetilde{R}'_-(h)+ \Re e \,\widetilde{R}^{\eps}_-(h), \\[12pt]
\ 0 &=& \dsp\eps^{-1} \Im m\,R_-^\eps(h) = h_1\,\Im m\, \frac{\partial R_-'}{\partial h_1}	(0,0)+ h_2\,\Im m\, \frac{\partial R_-'}{\partial h_2 }	(0,0) + \Im m \,\widetilde{R}'_-(h)+ \Im m \,\widetilde{R}^{\eps}_-(h).
\end{array}
\end{equation}
This is equivalent to have 
\begin{equation}\label{system1}
\bbM h=G^{\eps}(h)\qquad\mbox{ with }\qquad G^{\eps}(h):=\left(\begin{array}{c}
 -\Re e \,\widetilde{R}'_-(h)- \Re e \,\widetilde{R}^{\eps}_-(h)\\[4pt]
-\Im m \,\widetilde{R}'_-(h)- \Im m \,\widetilde{R}^{\eps}_-(h)
\end{array}\right)
\end{equation}
\begin{equation}\label{matrixM}
\mbox{ and }\qquad \bbM:=\left(\begin{array}{cc}
\Re e\, \cfrac{\partial R_-'}{\partial h_1}	(0,0) & \Re e\, \cfrac{\partial R_-'}{\partial h_2}	(0,0) \\[12pt]
\Im m\, \cfrac{\partial R_-'}{\partial h_1}	(0,0) & \Im m\, \cfrac{\partial R_-'}{\partial h_2}	(0,0)
\end{array}\right)\in\bbR^{2\times2}.
\end{equation}
To continue the procedure, we have to assume that the geometrical parameters $y_j^\pm, L_j, \ell_j, H_j^\pm$ are also such that the matrix $\bbM$ is invertible (again, see \S\ref{paragraphCalculus} for examples where this assumption is satisfied). Then $h$ solves (\ref{system1}) if and only if it is a solution to the problem
\begin{equation}\label{PbFixedPoint}
\begin{array}{|l}
\ \mbox{Find }h\in\bbR^2\mbox{ such that }\\[4pt]
\ h=\mathscr{G}^{\eps}(h)
\end{array}
\end{equation}
with $\mathscr{G}^{\eps}(h):=\bbM^{-1} G^{\eps}(h)$. Thus we obtain a fixed point equation. Let us prove that for a given $\varrho > 0$ small enough, $\mathscr{G}^{\eps}(\cdot)$ is a contraction in the closed ball $\overline{\bbB_\varrho} :=  \{ h \in \bbR^2 : 
|h| \leq \varrho \}$ for all $\eps\in(0;\eps_0(\varrho)]$. First, since $h\mapsto R'_-(h)$ is analytic, we have 
\begin{equation}\label{Analytic1}
|\widetilde{R}'_-(h)|\le c\,|h|^2\quad\mbox{ and }\quad|\widetilde{R}'_-(h^{(1)})-\widetilde{R}'_-(h^{(2)})|\le c\,|h^{(1)}+h^{(2)}|\,|h^{(1)}-h^{(2)}|,\quad \forall h,h^{(1)},h^{(2)}\in\overline{\bbB_\varrho}.
\end{equation}
We emphasize that $c$ in \ef{Analytic1} and in the estimates below is independent of $\varrho$ small enough. Moreover, we will prove that $h\mapsto R^{\eps}_-(h)$ is smooth. We will deduce that the remainder $h\mapsto \widetilde{R}^{\eps}_-(h)$ in (\ref{ExpansionEpsh}) is also smooth. More precisely, we will show (see Theorem \ref{FF} and \S\ref{sec5.3}) the estimates, for all $\delta\in(0;1/2)$,
\begin{equation}\label{Analytic2}
|\widetilde{R}^{\eps}_-(h)|\le c\,\eps^{1/2-\delta}\quad\mbox{ and }\quad|\widetilde{R}^{\eps}_-(h^{(1)})-\widetilde{R}^{\eps}_-(h^{(2)})|\le c\,\eps^{1/2-\delta}\,|h^{(1)}-h^{(2)}|,\qquad \forall h,h^{(1)},h^{(2)}\in\overline{\bbB_\varrho}.
\end{equation}
This allows us to write, using also the invertibility of $\bbM$, 
\[
|\mathscr{G}^{\eps}(h)| \le c\,|G^{\eps}(h)| \le c\,\big(|h|^2 
+ \eps^{1/2- \delta} \big)
\]
and guarantees that $\mathscr{G}^{\eps}(\cdot)$ maps $\overline{\bbB_\varrho}$ to $\overline{\bbB_\varrho}$ for $\varrho$ small enough and all $\eps\in(0;\eps_0(\varrho)]$. On the other hand, for $h^{(1)}$, $h^{(2)}\in\overline{\bbB_\varrho}$, we have 
\begin{equation}
\begin{array}{lcl}
|\mathscr{G}^{\eps}(h^{(1)})-\mathscr{G}^{\eps}(h^{(2)})| &\le &
c\,\big( \big| \widetilde{R}_-'(h^{(1)}) - \widetilde{R}_-'(h^{(2)}) \big| 
+ \big| \widetilde{R}_-^{\eps}(h^{(1)}) - \widetilde{R}_-^{\eps}(h^{(2)}) \big|\big) \\[4pt]
&\le & c\,\big(  |h^{(1)}| + |h^{(2)}|
+ \eps^{1/2- \delta} \big)\,|h^{(1)} -  h^{(2)} |,
\end{array}
\end{equation}
which ensures that $\mathscr{G}^{\eps}(\cdot)$ is a contraction for $\varrho$ small enough and all $\eps\in(0;\eps_0(\varrho)]$. Now the Banach contraction mapping principle proves the existence of a 
unique solution $h^{\mrm{sol}}$ to (\ref{PbFixedPoint}) in $\overline{\bbB_\varrho}$. Moreover, we have the estimate
\bea\label{R10}
|h^{\mrm{sol}}_1| + |h^{\mrm{sol}}_2| \leq c_ 0\,\eps^{1/2- \delta} .  
\eea
This leads us to the announced assertion on non-reflectability.

\BET\label{thNoRe}
Let $\sfk \geq 0$ and $\lambda > \lambda_\dagger$ be fixed. Assume that the geometrical parameters $y_j^\pm, L_j, \ell_j, H_j^\pm$, $j= 1, \ldots, J$ are such that the coefficient $R'_{-}$ in (\ref{MainCorR}) satisfies $R'_{-}=0$ and the matrix $\bbM$ in (\ref{matrixM}) is invertible. Then, there exist $\eps_0 > 0$ and $c_0 > 0$
such that, for all $\eps \in (0;\eps_0]$ the problem (\ref{PbFixedPoint}) has a solution $h^{\mrm{sol}} =(h^{\mrm{sol}}_1, h^{\mrm{sol}}_2)$ which satisfies the estimate \ef{R10}. Then we have $R^{\eps}_-(h^{\mrm{sol}})=0$, \textit{i.e.} in the channel $\Pi^{\eps}(h^{\mrm{sol}})$, the wave $w^+$ passes the family of tunnels without reflection. 
\ENT

\noindent Note that according to formula (\ref{MainCorT}), we have $\Im m\,T'(h^{\mrm{sol}})\le0$. And there holds $\Im m\,T'(h^{\mrm{sol}})<0$ as soon as the parameters are not such that $[\sfk=0\mbox{ and }w^{+}(y_j^-,-d)=w^{+}(y_j^+,-d),\ \forall j=1,\ldots,J]$. When $\Im m\,T'(h^{\mrm{sol}})<0$, at least for $\eps$ small enough, the transmission coefficient $T^\eps(h^{\mrm{sol}})$ has a negative imaginary part: 
$\Im m\,T^\eps(h^{\mrm{sol}})<0$. This means that after passing the tunnels
the wave $w^+$ certainly gets a phase shift, although it is of the order $\eps$ only.

\BER
To achieve non-reflectability we have only varied the positions of the tunnels, although the other geometric parameters \ef{R2} could evidently be varied as well. However, doing so, the smooth dependence of the reflection coefficient $R_-^\eps(h)$ on the perturbation parameter $h$ would require other arguments, which would possibly be more involved than the simple change of coordinates in  \ef{R6}. \ \ $\boxtimes$
\ENR

\subsection{The case of wells.}\label{sec5.4}
We mainly keep the  notation of Section \ref{sec1.1} but replaced the tunnels $\varpi_{j}^\eps$ by the wells 
\[
\varpi_{j-}^\eps := \{ x \in \varpi_j^\eps : s_j < 0 \},\quad\ j=1, \ldots, J. 
\]
Then we set $\Pi^\eps = \Pi \cup \varpi_{1 - }^\eps \cup \ldots
\varpi_{J - }^\eps$ (see Fig.\,\ref{fig2},\,b)) and consider the original water-wave problem \ef{3}--\ef{5} in the new channel $\Pi^\eps$. Then adapting the approach below, we find that the reflection and transmission coefficients for the scattering solution $u_-^{\eps}$ (see \ef{totalFields}) admits the asymptotic expansions 
\begin{equation}\label{MainAnsatzRefDemi}
R^{\eps}_-=0+\eps R^{'}_-+\dots;\qquad\qquad T^{\eps}=1+\eps T^{'}+\dots.
\end{equation}
\[
\mbox{with}\qquad\qquad\qquad R'_-=\cfrac{-i}{2N\sqrt{a^2-k^2}}\,\dsp\sum_{j=1}^J\int_{-\ell_j}^{0} H_j(s)\left(\bigg(\cfrac{dv^0_j}{ds}(s)\bigg)^2+\sfk^2(v^0_j(s))^2\right)\,ds
\]
\[
\mbox{and }\qquad\qquad\qquad T'=\cfrac{-i}{2N\sqrt{a^2-k^2}}\,\dsp\sum_{j=1}^J\int_{-\ell_j}^{0} H_j(s)\left(\bigg|\cfrac{dv^0_j}{ds}(s)\bigg|^2+\sfk^2|v^0_j(s)|^2\right)\,ds.
\]
Here again $N=(2a)^{-1}\,(e^{2ad}-e^{-2ad})+2d$ and $H_j(s) = H^+_j(s) + H^-_j(s) > 0$ corresponds to the rescaled thickness of the strip $\varpi_{j-}^\eps$. Moreover, the functions $v^0_j$  are defined as the solutions to the problems
\begin{equation}\label{MixedBVP}
\begin{array}{|lcl}
\ -\cfrac{d}{ds}\Big(H_j(s)\cfrac{dv^0_j}{ds} (s)\Big)+\sfk^2H_j(s)v^0_j(s)=0\qquad\mbox{for }s\in(-\ell_j;0)\\[12pt]
\hspace{3cm}v^0_j(- \ell_j)=w^{+}(y^-_j,-d),\hspace{1cm}\partial_sv^0_j(0)=0.
\end{array}
\end{equation}
Note that the last Neumann condition originates from  the no-penetration condition $\partial_\nu u^\eps(x) = 0$ on the end face $\{ x \in \varpi_j^\eps: s_j= 0\}$ of the curved well $\varpi_{j-}^\eps$. The mixed boundary-value problem \ef{MixedBVP} is still uniquely solvable. The justification of asymptotics in this case is completely similar to the case of tunnels below. As a consequence, the fine-tuning scheme of \S\ref{FineTunningParag} still works and provides examples of families of non-reflecting wells.

\section{Numerical experiments.}\label{SectionNumerics}
In this section, we implement numerically the approach leading to the Theorem \ref{thNoRe} above. 
\subsection{Preliminaries calculus.}\label{paragraphCalculus} The first step in the procedure presented in the previous section consists in finding geometrical parameters $y_j^\pm, L_j, \ell_j, H_j^\pm$, $j= 1, \ldots, J,$ such that the coefficient $R'_-$ in (\ref{MainCorR}) satisfies $R'_-=0$ and the matrix $\bbM$ in (\ref{matrixM}) is invertible. To proceed, we divide the analysis according to the value of $\sfk$.\\
\newline
$\star$ {\sc Case $\sfk=0$.} When $\sfk=0$, the solution of problem (\ref{Problem1D}) satisfies $H_j(s)d_sv^0_j (s)=\mrm{constant}$ for $s\in(-\ell_j;\ell_j)$. Using the boundary conditions, we deduce that 
\[
v^0_j(s)=w^{+}(y^-_j,-d)+\big(w^{+}(y^+_j,-d)-w^{+}(y^-_j,-d)\big)\Big(\dsp\int_{-\ell_j}^{\ell_j}\cfrac{1}{H_j(s)}\,ds\Big)^{-1}\Big(\dsp\int_{-\ell_j}^{s}\cfrac{1}{H_j(s)}\,ds
\Big). 
\]
Set $\mathcal{H}_j:=\Big(\dsp\int_{-\ell_j}^{\ell_j}\cfrac{1}{H_j(s)}\,ds\Big)^{-1}$. We have $d_sv^0_j(s) =(w^{+}(y^+_j,-d)-w^{+}(y^-_j,-d)) \mathcal{H}_j/H_j(s)$. Inserting the latter relation in (\ref{MainCorR}) and using that $w^{+}(y^\pm_j,-d)=2e^{iy^\pm_j\sqrt{a^2-k^2}}$, we obtain 
\begin{equation}\label{RefkNull}
R'_-=\cfrac{-i}{2Na}\,\dsp\sum_{j=1}^J\mathcal{H}_j(w^{+}(y^+_j,-d)-w^{+}(y^-_j,-d))^2=\cfrac{-2i}{Na}\,\dsp\sum_{j=1}^J\mathcal{H}_j(
e^{iy^+_ja}-e^{iy^-_ja})^2.
\end{equation}
Take $J=3$, $H_j=H$, $\ell_j=\ell$, $j=1,2,3$ so that $\mathcal{H}_j=\mathcal{H}$, $j=1,2,3$. In other words, we consider a situation with three similar tunnels. Then set 
\begin{equation}\label{ChoiceParam1}
y^-_1=-\cfrac{\pi}{3a};\qquad y^-_2=0;\qquad y^-_3=\cfrac{\pi}{3a};\quad\mbox{ and }\quad y^+_j=y^-_j+\cfrac{\eta}{a}\ \mbox{ with }\ \eta\in(0;\cfrac{\pi}{3}\,).
\end{equation}
In this case, according to (\ref{RefkNull}), we have 
\[
R'_-=\cfrac{-2i}{Na}\,\mathcal{H}\dsp\sum_{j=1}^3(e^{i\eta }-1)^2\,e^{2iy^-_ja}=\cfrac{-2i}{Na}\,\mathcal{H}\,(e^{i\eta }-1)^2(e^{-2i\pi/3}+1+e^{2i\pi/3})=0.
\]
Then, when we translate the position of the tunnels 1 and 2 respectively by $h_1$ and $h_2$, according to (\ref{RefkNull}), we find $R'_-(h_1,h_2)=-2i\mathcal{H}(Na)^{-1}\,(e^{i\eta}-1)^2\,(e^{2i(y^-_1+h_1)a}+e^{2i(y^-_2+h_2)a}+e^{2iy^-_3a})$. We deduce 
\[
\cfrac{\partial R'_-}{\partial h_1}(0,0)=\cfrac{4}{N}\,\mathcal{H}\,(e^{i\eta}-1)^2\,e^{2iy^-_1a}\qquad\mbox{ and }\qquad\cfrac{\partial R'_-}{\partial h_2}(0,0)=\cfrac{4}{N}\,\mathcal{H}\,(e^{i\eta}-1)^2\,e^{2iy^-_2a}.
\]
Thus the matrix $\bbM$ in (\ref{matrixM}) is invertible if and only if the  ratio $e^{2iy^-_2a}/e^{2iy^-_1a}=e^{2i\pi/3}$ is not a real number, which is indeed the case.\\
\newline
$\star$ {\sc Case $\sfk>0$.} When $\sfk>0$, to simplify the presentation, we assume that $H_j(s)=H^0_j$ does not depend on $s$ (the width of the tunnels is constant). Then the function $v^0_j$ satisfies 
\[
\begin{array}{|lcl}
\ -\cfrac{d^2v^0_j}{ds^2} (s)+\sfk^2v^0_j(s)=0\qquad\mbox{for }s\in(-\ell_j;\ell_j)\\[12pt]
\hspace{3cm}v^0_j(\pm \ell_j)=w^{+}(y^\pm_j,-d).
\end{array}
\]
We deduce that $v^0_j$ is given by
\[
v^0_j(s)=w^{+}(y^-_j,-d)\cosh(\sfk(s+\ell_j))+\big(w^{+}(y^+_j,-d)-w^{+}(y^-_j,-d)\cosh(2\sfk\ell_j)\big)\,\cfrac{\sinh(\sfk(s+\ell_j))}{\sinh(2\sfk\ell_j)}\ . 
\]
We can also write $v^0_j(s)=A_j\,e^{\sfk s}+B_j\,e^{-\sfk s}$ with 
\[
A_j=\cfrac{e^{\sfk\ell_j}e^{iy^+_j\sqrt{a^2-\sfk^2}}-e^{-\sfk\ell_j}e^{iy^-_j\sqrt{a^2-\sfk^2}}}{\sinh(2\sfk\ell_j)}\quad\mbox{ and }\quad B_j=\cfrac{-e^{-\sfk\ell_j}e^{iy^+_j\sqrt{a^2-\sfk^2}}+e^{\sfk\ell_j}e^{iy^-_j\sqrt{a^2-\sfk^2}}}{\sinh(2\sfk\ell_j)}.
\]
Then we have
\[
\int_{-\ell_j}^{\ell_j} \bigg(\cfrac{dv^0_j}{ds}(s)
\bigg)^2+\sfk^2(v^0_j(s))^2\,ds=2\sfk(A_j^2+B_j^2)\sinh(2\sfk\ell_j).
\]
Using in particular that $H_j=2\ell_j\mathcal{H}_j$, we deduce
\begin{equation}\label{SNonZerok}
\begin{array}{lcl}
R'_-&=&\cfrac{-i}{2N\sqrt{a^2-\sfk^2}}\,\dsp\sum_{j=1}^J\int_{-\ell_j}^{\ell_j} H_j\bigg(\Big(\cfrac{dv^0_j}{ds}(s)\Big)^2+\sfk^2(v^0_j(s))^2\bigg)\,ds\\[16pt]
&=&\cfrac{-i\sfk}{N\sqrt{a^2-\sfk^2}}\,\dsp\sum_{j=1}^JH_j(A_j^2+B_j^2)\sinh(2\sfk\ell_j)\\[16pt]
&=&\cfrac{-2i}{N\sqrt{a^2-\sfk^2}}\,\dsp\sum_{j=1}^J\mathcal{H}_j\Big(\mathcal{A}_j\,(e^{2iy^-_j\sqrt{a^2-\sfk^2}}+e^{2iy^+_j\sqrt{a^2-\sfk^2}})-2\mathcal{B}_j\,e^{i(y^-_j+y^+_j)\sqrt{a^2-\sfk^2}}\,\Big)
\end{array}
\end{equation}
with
\[
\mathcal{A}_j =2\sfk\ell_j\,\cfrac{e^{2\sfk\ell_j}+e^{-2\sfk\ell_j}}{e^{2\sfk\ell_j}-e^{-2\sfk\ell_j}}\qquad\mbox{ and }\qquad \mathcal{B}_j =\cfrac{4\sfk\ell_j}{e^{2\sfk\ell_j}-e^{-2\sfk\ell_j}}.
\]
Note that one can verify that taking the limit $\sfk\to0$ in (\ref{SNonZerok}), we get back the relation (\ref{RefkNull}) obtained for $\sfk=0$. In order to have $R'_-=0$, we take $J=3$, $H_j=H$, $\ell_j=\ell$, $j=1,2,3$ so that $\mathcal{A}_j= \mathcal{A},\ \mathcal{B}_j= \mathcal{B},\ \mathcal{H}_j=\mathcal{H},\,j=1,2,3$. In other words, again we consider three similar tunnels. Then we set 
\begin{equation}\label{ChoicePositionkPos}
\begin{array}{lll}
\dsp y_1^- = - \frac{19}{12} \frac{\pi}{\sqrt{a^2 -\sfk^2 }}, \qquad &\dsp y_2^- = \frac34 \frac{\pi}{\sqrt{a^2 -\sfk^2 }}, \qquad & \dsp y_3^- =  \frac{25}{12} \frac{\pi}{\sqrt{a^2 -\sfk^2 }}, \\ [12pt]
\dsp y_1^+ = - \frac{13}{12} \frac{\pi}{\sqrt{a^2 -\sfk^2 }}, \qquad&
\dsp y_2^+ = \frac54 \frac{\pi}{\sqrt{a^2 -\sfk^2 }}, \qquad &
\dsp y_3^+  =  \frac{31}{12} \frac{\pi}{\sqrt{a^2 -\sfk^2 }}.  
\end{array}
\end{equation}
Using that $y_j^+-y_j^-=\pi/(2\sqrt{a^2 -\sfk^2 })$, we find 
\begin{equation}\label{SNonZerok}
\begin{array}{lcl}
R'_-&=&\cfrac{-2i}{N\sqrt{a^2-\sfk^2}}\,\mathcal{H}\Big(-2\mathcal{B}(e^{-2i\pi/3}+0+e^{2i\pi/3})+\dsp\sum_{j=1}^3\mathcal{A}\,e^{2iy^-_j\sqrt{a^2-\sfk^2}}\,(1-1)\,\Big)=0.
\end{array}
\end{equation}
Now translating the position of the tunnels 1 and 2 respectively by $h_1$ and $h_2$, according to (\ref{RefkNull}), we find 
\[
\begin{array}{lcl}
R'_-(h_1,h_2)=\cfrac{-2i}{N\sqrt{a^2-\sfk^2}}\,\mathcal{H}\dsp\sum_{j=1}^3\Big(\mathcal{A}\,(e^{2i(y^-_j+h_j)\sqrt{a^2-\sfk^2}}+e^{2i(y^+_j+h_j)\sqrt{a^2-\sfk^2}})\\[14pt]
\hspace{8cm}-2\mathcal{B}\,e^{i((y^-_j+h_j)+(y^+_j+h_j))\sqrt{a^2-\sfk^2}}\,\Big)
\end{array}
\]
with $h_3=0$. Then for $j=1,2$, we obtain
\[
\cfrac{\partial R'_-}{\partial h_j}(h_1,h_2)=\cfrac{4}{N}\,\dsp\mathcal{H}\Big(\mathcal{A}\,(e^{2i(y^-_j+h_j)\sqrt{a^2-\sfk^2}}+e^{2i(y^+_j+h_j)\sqrt{a^2-\sfk^2}})-2\mathcal{B}\,e^{i((y^-_j+h_j)+(y^+_j+h_j))\sqrt{a^2-\sfk^2}}\,\Big)
\]
and so
\[
\cfrac{\partial R'_-}{\partial h_j}(0,0)=\cfrac{4}{N}\,\dsp\mathcal{H}\Big(\mathcal{A}\,(e^{2iy^-_j\sqrt{a^2-\sfk^2}}+e^{2iy^+_j\sqrt{a^2-\sfk^2}})-2\mathcal{B}\,e^{i(y^-_j+y^+_j)\sqrt{a^2-\sfk^2}}\,\Big).
\]
Thus for our choice (\ref{ChoicePositionkPos}) of the positions of the tunnels, we obtain
\[
\cfrac{\partial R'_-}{\partial h_1}(0,0)=-\cfrac{8}{N}\,\mathcal{H}\,\mathcal{B}\,e^{-2i\pi/3} \qquad\mbox{ and }\qquad \cfrac{\partial S'_-}{\partial h_2}(0,0)=-\cfrac{8}{N}\,\mathcal{H}\,\mathcal{B}\,.
\] 
As a consequence, the matrix $\bbM$ is indeed invertible.

\BER\label{rem5R}
If $J=2$, which corresponds to the case of only two tunnels, the matrix
$\bbM$ is always non invertible because $|R^\eps_-|$ is invariant with respect to the coordinate change $(y,z) \mapsto (y - 
h^0,z)$. In this case we would only have the real parameter $h_2-h_1$
instead of the couple $(h_1,h_2)$ which is not enough to cancel one complex coefficient.  \ \ $\boxtimes$
\ENR
\noindent Above we presented two examples. Of course, the list can be enlarged readily.
\subsection{Numerical procedure.} Numerically, we solve the fixed point problem (\ref{PbFixedPoint}) using an iterative procedure. More precisely, we start from some arbitrary $h^{0}:=(h_1^0,h_2^0)^{\top}\in\bbR^2$ whose norm is not ``too large'' and then, for a given $h^n\in\bbR^2$, we set 
\begin{equation}\label{FixedPointNum}
h^{n+1}:=\mathscr{G}^{\eps}(h^n)\qquad \Leftrightarrow\qquad \bbM h^{n+1}=G^{\eps}(h^n)=\left(\begin{array}{c}
 -\Re e \,\widetilde{R}'_-(h^n)- \Re e \,\widetilde{R}^{\eps}_-(h^n)\\[4pt]
-\Im m \,\widetilde{R}'_-(h^n)- \Im m \,\widetilde{R}^{\eps}_-(h^n)
\end{array}\right).
\end{equation}
Let us explain how to compute the right hand side in (\ref{FixedPointNum}). Using (\ref{Expansionh}) in (\ref{ExpansionEpsh}), we get
\[
R^{\eps}_-(h^n)=\eps \left(h^n_1\cfrac{\partial R_-'}{\partial h_1}(0,0)+h^n_2\cfrac{\partial R_-'}{\partial h_2}(0,0)+\widetilde{R}_-'(h^n)+\widetilde{R}^{\eps}_-(h^n) \right).
\]
Extracting the real and imaginary parts, we deduce from the definition of $G^{\eps}(h^n)$ (see (\ref{FixedPointNum})) that there holds
\[
G^{\eps}(h^n)=\bbM h^n-\eps^{-1}\left(\begin{array}{c}
\Re e\,R^{\eps}_-(h^n)\\
\Im m\,R^{\eps}_-(h^n)
\end{array}\right).
\]
Using the latter relation in (\ref{FixedPointNum}), we see that $h^{n+1}$ is related to $h^n$ by the simple formula
\begin{equation}\label{MainEquationNum}
h^{n+1}:=h^{n}-\eps^{-1}\bbM^{-1}\left(\begin{array}{c}
\Re e\,R^{\eps}_-(h^n)\\
\Im m\,R^{\eps}_-(h^n)
\end{array}\right).
\end{equation}
In (\ref{MainEquationNum}), the coefficient $R^{\eps}_-(h^n)$ is computed at each step $n\ge0$ solving the scattering problem 
\begin{equation}\label{PbChampTotalNumericalSteps}
\begin{array}{|rcll}
\multicolumn{4}{|l}{\mbox{Find }u^{\eps}_-(h^n)\mbox{ such that }u^{\eps}_-(h^n)-w^+\mbox{ is outgoing and } }\\[3pt]
- \Delta u^{\eps}_-(h^n) + {\sf k}^2 u^{\eps}_-(h^n) &=&  0 &\mbox{ in }\Pi^\eps(h^n)\\[3pt]
\partial_\nu u^{\eps}_-(h^n) &=& 0 &\mbox{ on }\partial \Pi^\eps(h^n)
\setminus \Gamma_0\\[3pt]
\partial_\nu u^{\eps}_-(h^n) &=& \lambda u^{\eps}_-(h^n) &\mbox{ on }\Gamma_0.\end{array}
\end{equation}
Note that the geometry of the channel $\Pi^{\eps}(h^n)$ depends on the step $n\ge0$. Then according to the representation (\ref{totalFields}), the coefficient $R^{\eps}_-(h^n)$ is given by
\begin{equation}\label{DefTermeChampLointain}
R^{\eps}_-(h^n)  = N^{-1}\dsp\int_{\{-2\ell_0\}\times(-d;0)} (u^{\eps}_-(h^n)-w^+)\,w^+\,d\sigma,
\end{equation}
where $N=(2a)^{-1}\,(e^{2ad}-e^{-2ad})+2d$ (see after (\ref{MainCorT})). At each step $n\ge0$, we approximate the solution of problem (\ref{PbChampTotalNumericalSteps}) with a P2 finite element method in $\Pi^{\eps}_{11}(h^n):=\{x=(y,z)\in\Pi^{\eps}(h^n)\,|\,|y|<11\}$. We emphasize in particular that at each step, it is necessary to mesh the domain. At $y=\pm 11$, a truncated Dirichlet-to-Neumann map with 10 terms serves as a transparent boundary condition. Computations are implemented with \textit{FreeFem++}\footnote{\textit{FreeFem++}, \url{http://www.freefem.org/ff++/}.} while  results are displayed with  \textit{Paraview}\footnote{\textit{Paraview}, \url{http://www.paraview.org/}.}. 

\subsection{Results.} 
\begin{figure}[!ht]
\centering
\includegraphics[width=14cm]{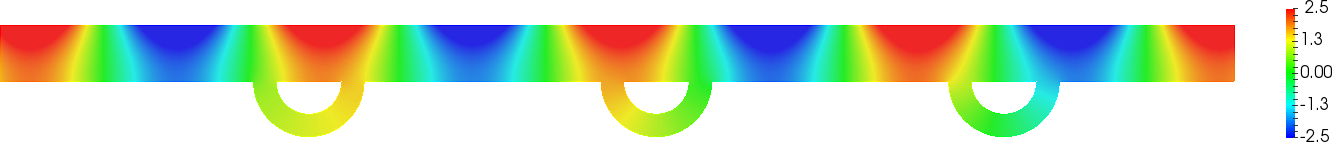}\\[12pt]\includegraphics[width=14cm]{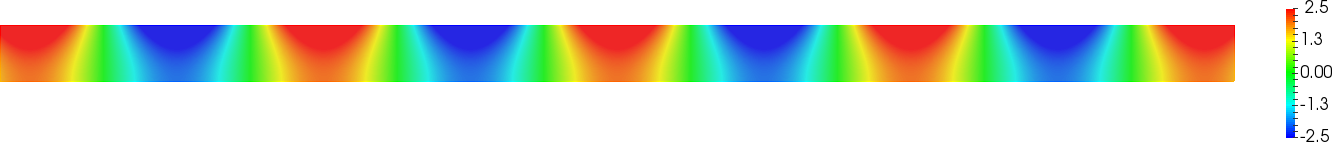}\\[12pt]
\includegraphics[width=14cm]{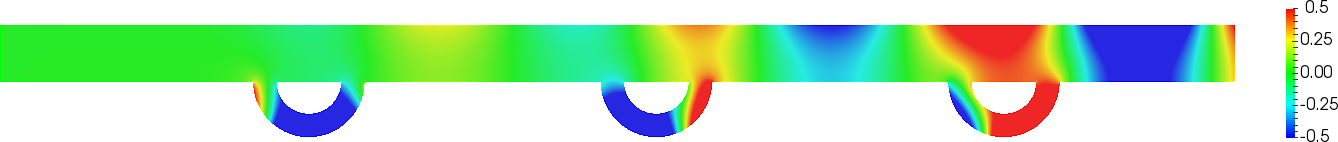}
\caption{Real parts of $u^{\eps}_-(h^n)$ (top), $w^+$ (middle) and $u^{\eps}_-(h^n)-w^+$ (bottom) after 15 iterations. Here $\sfk=0$ and $\lambda=1$. We obtain $R^{\eps}_-\approx(0.17+1.87\mrm{i})\,10^{-6}$ and $h^{\mrm{sol}}\approx(-0.18,-0.09)^{\top}$. 
\label{figResult1}}
\end{figure}

\begin{figure}[!ht]
\centering
\includegraphics[width=14cm]{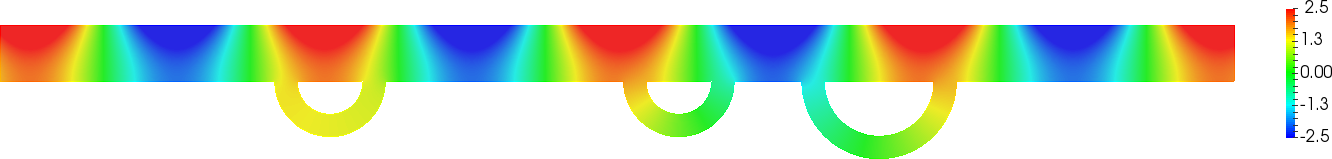}\\[12pt]\includegraphics[width=14cm]{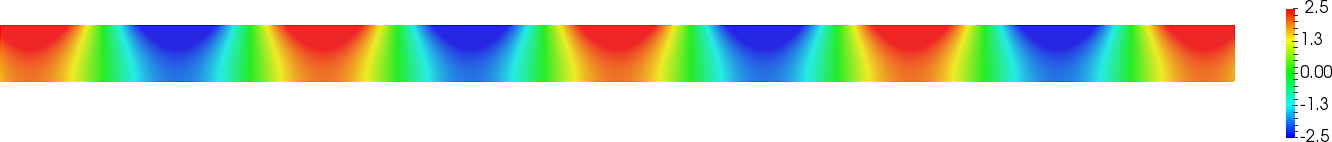}\\[12pt]
\includegraphics[width=14cm]{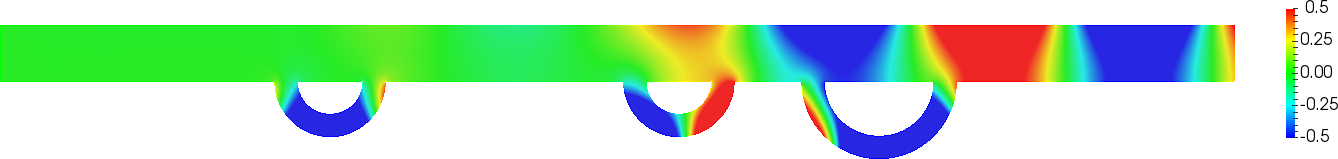}
\caption{Real parts of $u^{\eps}_-(h^n)$ (top), $w^+$ (middle) and $u^{\eps}_-(h^n)-w^+$ (bottom) after 15 iterations. Here $\sfk=0$ and $\lambda=1$. We obtain $R^{\eps}_-\approx(2.65-0.04\mrm{i})\,10^{-7}$ and $h^{\mrm{sol}}\approx(0.2,0.31)^{\top}$. 
\label{figResult2}}
\end{figure}

\begin{figure}[!ht]
\centering
\includegraphics[width=14cm]{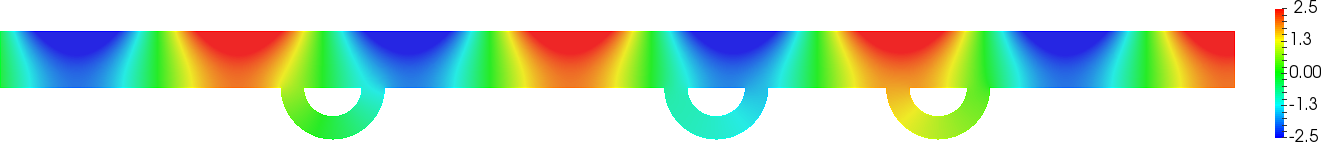}\\[12pt]\includegraphics[width=14cm]{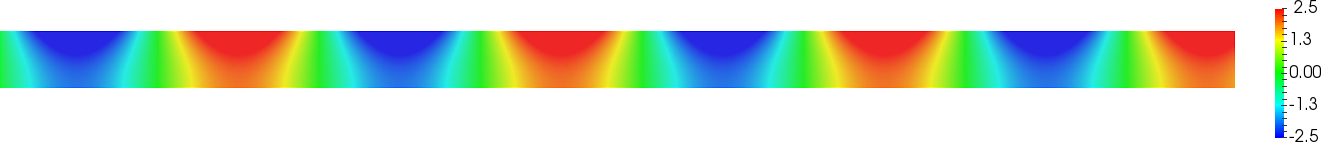}\\[12pt]
\includegraphics[width=14cm]{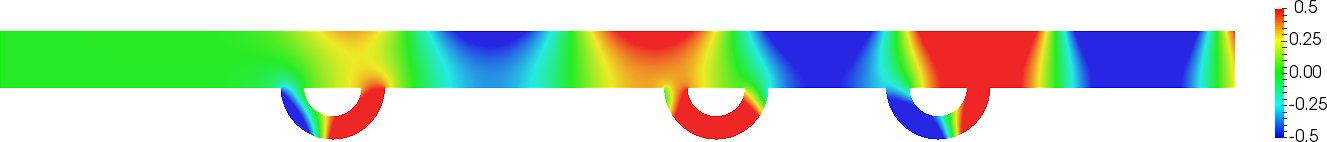}
\caption{Real parts of $u^{\eps}_-(h^n)$ (top), $w^+$ (middle) and $u^{\eps}_-(h^n)-w^+$ (bottom) after 15 iterations. Here $\sfk=0.5$ and $\lambda=1$. We obtain $R^{\eps}_-\approx(1.21+0.6\mrm{i})\,10^{-6}$ and $h^{\mrm{sol}}\approx(-0.23,-0.11)^{\top}$. 
\label{figResult3}}
\end{figure}

The results we obtain are displayed in Figures \ref{figResult1}--\ref{figResult3}. We start with $h^0=(0,0)^{\top}$ and set $\eps=0.2$. For the Figure \ref{figResult1}, we take $\sfk=0$, $\lambda=1$ and we fix the geometrical parameters as in (\ref{ChoiceParam1}) with $y^-_1=-\pi/(3a)-2\pi/a$, $y^-_2=0$, $y^-_3=\pi/(3a)+2\pi/a$ and $y^+_j=y^-_j+\pi/2$. For the Figure \ref{figResult2}, we set the parameters as in Figure \ref{figResult2} except for the definition of the third tunnel. Here  $\varpi_3^\eps$ is constructed from the half circle passing through $P_3^{\pm} = (y_3^{\pm},-d)$ with $y_3^{-}=\pi/(3a)+\pi/a$ and $y_3^{+}=y_3^{-}+3\pi/4$ (we enlarge the radius of the third tunnel to break the symmetry). We emphasize that in this setting we do not have $R'_-(0,0)=0$. However, this is not a problem and one can prove that the sequence $(h^n)$ constructed via the recursive relation (\ref{MainEquationNum}) converges to some $h^{\mrm{sol}}\in\bbR^2$ such that $R^{\eps}_-(h^{\mrm{sol}})=0$. Finally, for the Figure \ref{figResult3}, we take $\sfk=0.5$, $\lambda=1$ and the geometrical parameters as in (\ref{ChoicePositionkPos}). For each setting, we represent the real parts of $u^{\eps}_-(h^n)$ (top), $w^+$ (middle) and $u^{\eps}_-(h^n)-w^+$ (bottom) after 15 iterations ($n=15$). As expected, we observe that the amplitude of the scattered field $u^{\eps}_-(h^n)-w^+$ is very small in the incident direction (no reflection). Interestingly, the fixed point procedure converges though the parameter $\eps$ (the width of the tunnels) is not very small (here $\eps=0.2$). As predicted by the theory, one can also observe a small phase shift between the incident and the transmitted fields.

\section{Derivation of the asymptotic expansions}\label{SectionAsymptotics}
In this section, we compute an asymptotic expansion of the function $u^{\eps}_-$ defined in (\ref{totalFields}) with respect to $\eps$. In the straight channel $\Pi$, at least ``far'' from the junctions with the tunnels, we make the ansatz
\begin{equation}\label{MainAsymptoPi}
u^{\eps}_-=u^0+\eps u'+\ldots.
\end{equation}
Here and in what follows, the dots stand for inessential higher-order terms. Our goal is to identify the functions $u^0$, $u'$. This will allow us to compute the terms in the ansatz for the reflection and transmission coefficients
\[
R^{\eps}_-=R^0_-+R'_-+\ldots\ ;\qquad\qquad T^{\eps}=T^0+T'+\ldots.
\]
The justification of the expansions with the proof of error estimates will be given in Section \ref{SectionJustif}.

\subsection{Main asymptotic terms.} \label{sec3.1}
In $\Pi$, as a first approximation of $u^{\eps}_-$, it is natural to take $w^+$ because when $\eps$ tends to zero, the tunnels disappear and then the incident wave does not suffer from scattering. Therefore we set $u^0=w^+$.\\
\newline
Now let us focus our attention on the approximation of $u^{\eps}_-$ in the thin tunnels $\varpi_j^\eps$. To simplify the notation, we omit the index $j$. As usual, we stretch the transversal section of the tunnel considering the change of variables 
\bea
n \mapsto \zeta = \eps^{-1 }n \in  \Upsilon(s) := (\,-H^-(s);H^+(s) \,) \label{11}
\eea
while keeping the longitudinal coordinate $s$ unchanged.  The Laplace operator in the local coordinates 
\bea\label{12}
\Delta_x = A(n,s)^{-1} \partial_n (A(n,s) \partial_n\ ) +
A(n,s)^{-1} \partial_s (A(n,s)^{-1}  \partial_s \  ) 
\eea
admits the decomposition 
\bea
\Delta_x = \eps^{-2} \partial_\zeta^2 + \eps^{-1} \kappa(s)
\partial_\zeta + \partial_s^2 - \kappa(s)^2 \zeta \partial_\zeta + \ldots .
\label{13}
\eea
In (\ref{12}), $\partial_n = \partial / \partial n$ and so on, $A(n,s) = 1+ \kappa(s) n$ and $\kappa(s)$ is the curvature of $L$ at the point $s \in (-\ell;\ell)$. The unit outward normal vector $\nu_\pm^\eps$ on the lateral sides of the tunnel $\varpi^{\eps}$ is as follows:
\beas
\nu_\pm^\eps(s) = \big(1 + \eps^2 A(\pm \eps
H^\pm (s) , s )^{-2} |\partial_s H^\pm (s) |^2 \big)^{-1/2}
\bigg( \pm 1, - \cfrac{\eps \partial_sH^\pm (s)}{A(\pm \eps
H^\pm (s) , s )} \bigg). 
\eeas
Notice that this vector is written in the local coordinate system $(n,s)$ and therefore
\bea
\partial_{\nu_\pm^\eps} = \nu_\pm^\eps (s) \cdot 
\big(\partial_n , A(n,s)^{-1} \partial_s\big) 
= \pm \eps^{-1} \partial_\zeta- \eps \partial_s 
H^\pm(s) \partial_s + \ldots . \label{14}
\eea
Following the standard dimension reduction procedure in  thin domains, see e.g. \cite[Chap.\,15]{MaNaPl}, in $\varpi^\eps $, we consider the ansatz
\bea\label{15}
u^\eps(x) = v^0 (s) + \eps v'(\zeta,s) + \eps^2 
v''(\zeta,s) + \ldots 
\eea
where the functions $v^0$, $v'$ and $v''$ have to be determined. Inserting \ef{15} into the initial problem \ef{3}--\ef{5} restricted to $\varpi^\eps$, using \ef{13}, \ef{14} and collecting the terms of order $\eps^{-1}$, we arrive at 
\[
\begin{array}{|lcl}
\ - \partial_\zeta^2 v'(\zeta,s) &=& 0, \quad\zeta \in \Upsilon(s), \\
\ \pm \partial_\zeta v'( \pm H^\pm(s),s) &=&0 .
\end{array}
\]
Since the problem is homogeneous, we have to set $v'(\zeta ,s) = v'(s)$. Now collecting the terms at order $\eps^0=1$ when inserting \ef{15} in \ef{3}--\ef{5} restricted to $\varpi^\eps$, we obtain
\begin{equation}\label{16N}
\begin{array}{|lcl}
\ - \partial_\zeta^2 v''(\zeta,s) &=& \partial_s^2 v^0(s)
 - \sfk^2 v^0(s), \quad\zeta \in \Upsilon(s), \\
\ \pm \partial_\zeta v''( \pm H^\pm(s),s) &=& \partial_s H^\pm(s) \partial_s v^0 (s).
\end{array}
\end{equation}
The compatibility condition for \ef{16N} writes
\begin{equation}\label{eqn1Dint}
\int_{\Upsilon(s) } \partial_s^2 v^0(s)
 - \sfk^2 v^0(s)\,d\zeta + \partial_s H^+(s) \partial_s v^0 (s) + \partial_s H^-(s) \partial_s v^0 (s) = 0.
\end{equation}
Since the length of $\Upsilon(s)$ is equal to $H(s) = H^+(s) + H^-(s) > 0$ (the rescaled thickness of the curved strip $\varpi^\eps$), (\ref{eqn1Dint}) turns into the ordinary differential equation
\bea\label{17}
\ -\cfrac{d}{ds}\Big(H(s)\cfrac{dv^0}{ds} (s)\Big)+\sfk^2H(s)v^0(s)=0\qquad\mbox{for }s\in(-\ell;\ell).  
\eea
In order to close problem \ef{17}, we have to impose boundary conditions. To proceed, we match the value of $v^0$ with the one of $u^{\eps}_-|_{\Pi}$ at order $\eps^0=1$ at the junction points $P^{\pm}$. This gives us
\bea\label{18}
v^0(\pm \ell) = w^{+}(y^\pm,-d).
\eea
Thus \ef{17}-\ef{18} form the resultant problem on the arc $L$ introduced in \ef{Problem1D}. 
\subsection{Correction terms. } \label{paragraphCorrection}

\begin{figure}[!ht]
\centering
\includegraphics[width=10cm]{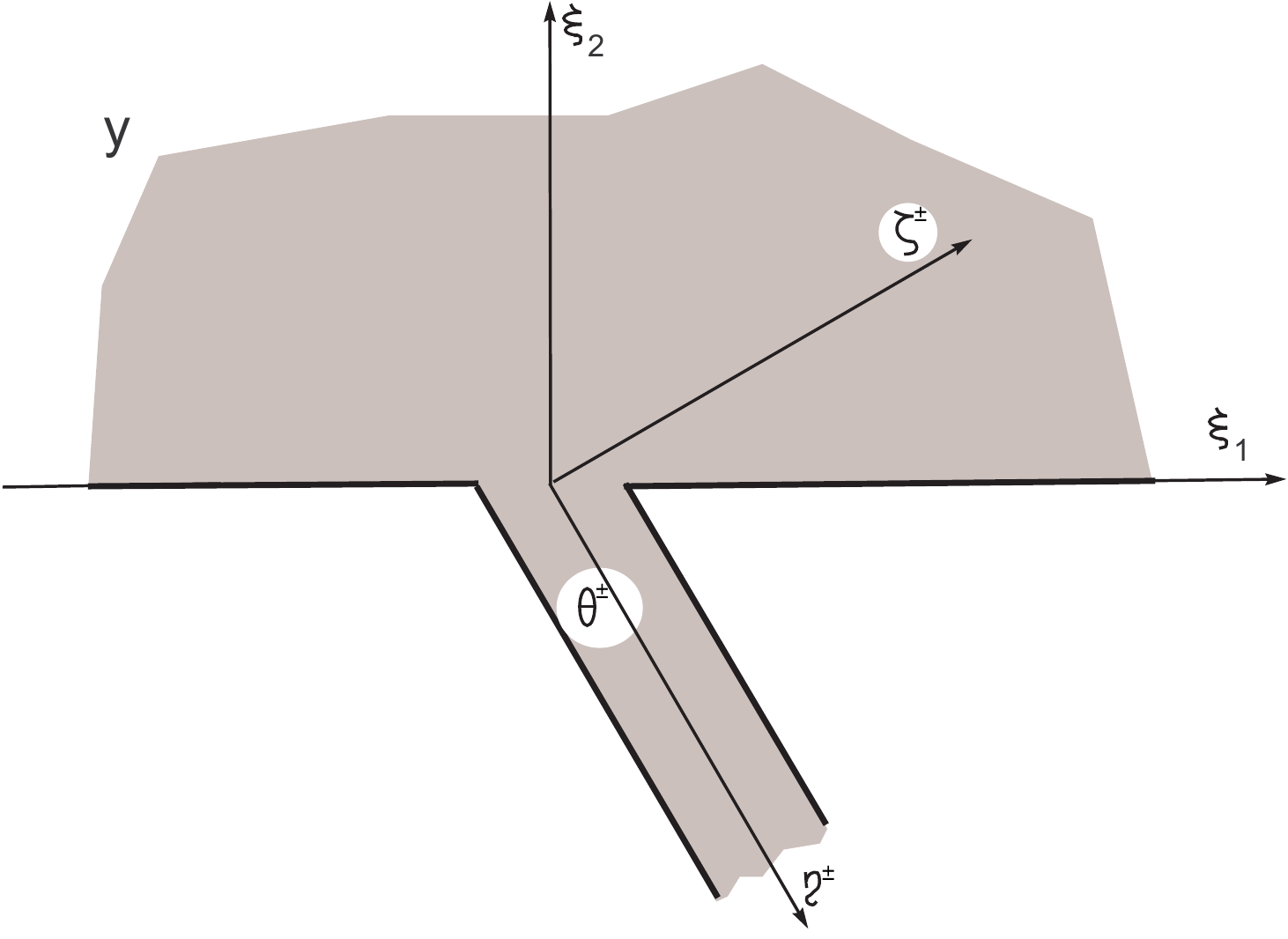}
\caption{Junction of the massive body and of the thin inclined tunnel.\label{fig3}}
\end{figure}

For the moment, we have formally derived an expansion at order $\eps^0=1$ of $u^{\eps}_-$ in $\Pi^{\eps}$. Now, we wish to get a better  approximation at order $\eps$. In particular, we want to identify the term $u'$ in the ansatz $u^{\eps}_-=w^++\eps u'+\ldots$ in $\Pi$ (see  (\ref{MainAsymptoPi})). Inserting this expansion in \ef{3}-\ef{5} and making $\eps\to0$, we find that $u'$ must satisfy 
\begin{equation}\label{E7}
\begin{array}{|rlcl}
\  -\Delta u'+\sfk u' &=& 0 &\qquad\mbox{in }\Pi, \\[3pt]
\ \partial_\nu  u'  &= &0& \qquad  \mbox{on }\Gamma_d\setminus {\textstyle \bigcup\limits_{\pm}} 
{\textstyle\bigcup\limits_{j=1}^J} P_j^\pm  \\[0pt]
\multicolumn{3}{|r}{\ \partial_\nu u' = \lambda u'} &\qquad\mbox{on }\Gamma_0.
\end{array}
\end{equation}
In order to define completely (uniquely) $u'$, we have to impose conditions at the junction points $P^\pm_j$. To proceed, we will employ the method of matched asymptotic expansions, cf. \cite{VD,Ilin}, which consists in matching the behaviour of $u^{\eps}_-$ in a neighbourhood of $P^\pm_j$ with the behaviour of a inner expansion of $u^{\eps}_-$ at infinity. Note that to obtain a non zero corrector $u'$, we have to allow singular behaviours at the $P^\pm_j$. In the following, to simplify notation, again we omit the index $j$.\\
\newline
In the junction zone near the endpoints $P^\pm$ of the curve $L$, we anticipate the existence of a boundary layer phenomenon. To capture it, in this region the asymptotic behaviour of $u^{\eps}_-$ will be described in the stretched coordinates 
\bea  \label{24}
\xi^\pm = (\xi_1^\pm, \xi_2^\pm) = \eps^{-1}( y- y^\pm,
z+d )=\eps^{-1} (x- P_j^\pm).
\eea
To simplify the proofs, especially the justification procedure in Section
\ref{SectionJustif}, cf. Remark \ref{rem2.1},  we assume that the profile functions $H^\pm $ in \ef{1} 
are constants and that the curve $L$ is straight for $s \in (-\ell-\delta; -\ell + \rho)$
and $s \in (\ell- \rho; \ell + \delta)$ for some $\rho > 0$. Thus, we see that changing the coordinate $x \mapsto \xi^{\pm}$ and setting formally  $\eps =0$ transform the domain $\Pi^{\eps}$ into the union of the upper half-plane $\bbR_+^2 := \{ \xi = (\xi_1, \xi_2) : \xi_2 > 0 \}$
and the rotated strip
\bea\label{mi1}
\theta^\pm := \{ \xi:(\eta^\pm,\zeta^\pm):=(\xi_1\cos\alpha^{\pm}-\xi_2\sin\alpha^{\pm},\xi_1\sin\alpha^{\pm}+\xi_2\cos\alpha^{\pm})\in \bbR\times\Upsilon(\pm \ell) \}
\eea
(see Figure \ref{fig3}). In the vicinity of the point $P^{\pm}$, we introduce the \textit{inner} expansion of $u^{\eps}_-$:
\bea\label{E1}
u^\eps_-(x) = \cU^{0 \pm} (\xi^\pm) + 
\eps \cU'^{\pm} (\xi^\pm) + \ldots . 
\eea
Since $\Delta_x - \sfk^2 = \eps^{-2} (\Delta_{\xi^\pm} - \eps^2 \sfk)$, the Laplacian is the main asymptotic part of the Helmholtz 
operator in the variables $\xi^{\pm}$. Furthermore,
on the boundary $\partial \Pi^\eps$ near 
$P^\pm$, the normal derivative $\partial_{\nu(x)}$ is nothing but 
$\eps^{-1} \partial_{\nu(\xi^\pm)}$, where $\nu(\xi^\pm)$ is the
outward normal on the boundary of the junction $\Theta^\pm := \bbR^2_+ \cup \theta^\pm$. Therefore, we deduce that the functions $\cU^{0 \pm}$, $\cU'^{\pm}$ must satisfy 
\bea\label{My27}
\begin{array}{|rlcl}
\  -\Delta_\xi \cU^{0\pm} &=& 0& \quad\mbox{ in } \Theta^{\pm} \\[3pt]
\ \partial_{\nu} \cU^{0\pm} &=& 0&\quad\mbox{ on } \partial\Theta^{\pm} 
\end{array}\qquad\qquad \begin{array}{|rlcl}
\  -\Delta_\xi \cU'^{\pm} &=& 0& \quad\mbox{ in } \Theta^{\pm} \\[3pt]
\ \partial_{\nu} \cU'^{\pm} &=& 0&\quad\mbox{ on } \partial\Theta^{\pm}. 
\end{array}
\eea
According to the conditions one imposes at infinity in $\Theta^{\pm}$, Problems \ef{My27} can admit non zero solutions. Clearly $\cY^0(\xi) = 1$ is one of them. One can also show (see Proposition \ref{PropositionExistenceTerm1} below) that there are some functions $\cY^{1\pm}$, $\cY^{2\pm}\in H^1_{\mrm{loc}}(\overline{\Theta^{\pm}})$ satisfying \ef{My27} such that when $|\xi| \to + \infty$, 
\begin{equation}\label{30}
 \cY^{1\pm} (\xi) = \left\{
\begin{array}{ll}
 \dsp\frac{1}{\pi} \ln \frac{1}{|\xi|} + O(|\xi|^{-1}) &\qquad \mbox{in }\bbR_+^2 \\[7pt]
\dsp\frac{\eta^{\pm}}{H(\pm\ell)}  + K^{1\pm } + O(e^{- \pi \eta^{\pm}/ H(\pm\ell)} )&\qquad\mbox{in }\theta^{\pm}\setminus\bbR_+^2 
\end{array}
\right.
\end{equation}
and
\begin{equation}\label{30N}
 \cY^{2\pm} (\xi) = \left\{
\begin{array}{ll}
\dsp\xi_1 +O(|\xi|^{-1})& \qquad\mbox{in }\bbR_+^2 \\[7pt]
\dsp K^{2\pm } + O(e^{- \pi \eta^{\pm}/ H(\pm\ell)})
&\qquad\mbox{in }\theta^{\pm}\setminus\bbR_+^2. 
\end{array}
\right.
\end{equation}
 Here, $K^{1\pm }$, $K^{2\pm }$ are constants depending on the width and the tilt angle of the strip $\theta^{\pm}$. Moreover $H^1_{\mrm{loc}}(\overline{\Theta^{\pm}})$ denotes the set of functions whose $H^1$ norm is finite in any bounded subdomain of $\overline{\Theta^{\pm}}$. We will look for $U^{0\pm}$, $U^{1\pm}$ as combinations of the $\cY^{0}$, $\cY^{1\pm}$, $\cY^{2\pm}$. One could add other solutions of \ef{My27} in the expansions but this is not needed for the calculus of the first two terms.\\
\newline
Now we match the behaviours of the two outer expansions of $u^{\eps}_-$ in $\Pi$ and in $\varpi^{\eps}$ with the behaviour of the inner expansion of $u^{\eps}_-$ in an intermediate region where $|x- P^{\pm}|\to0$ and $|\xi^{\pm}|\to+\infty$. Using the stretched coordinates \ef{24}, a Taylor expansion at the point $P^{\pm}$ yields
\begin{equation} \label{41}
\begin{array}{lcl}
w^+ (y,z) &=& w^+(y^\pm,-d) +
(y-y^\pm) \partial_y w^+ (y^\pm,-d)+
(z+d) \partial_z w^+(y^\pm,-d) + O((r^\pm)^2) \\[3pt]
 &=&w^+(y^\pm,-d) + \eps\xi_1^{\pm} \partial_y w^+(y^\pm,-d)
+ O(\eps^2).
\end{array}
\end{equation} 
Here, we took into account that $\partial_z w^+(y^\pm,-d)=0$
and set $r^\pm := |x- P^\pm| = ( |y- y^\pm|^2 +|z + d|^2 )^{1/2}$. On the other hand, in the vicinity of $P^{\pm}$ we also have
\begin{equation} \label{44}
\begin{array}{lcl}
v^0(s) &=& v^0(\pm \ell ) + s \partial_s v^0 (\pm \ell)
+ O(s^2) \\[3pt]
 &=&v^0(\pm \ell ) \mp \eps \eta^\pm
\partial_s v^0 (\pm \ell) + O(\eps^2) . 
\end{array}
\end{equation}
Identifying the powers in $\eps^0=1$, first, we obtain 
\bea\label{E2}
\cU^{0 \pm} (\xi^\pm) = w^+(y^\pm,-d)\cY^{0}(\xi^\pm)=w^+(y^\pm,-d) = v^0(\pm \ell_j). 
\eea
In the inner expansion $u^\eps_-(x) = \cU^{0 \pm} (\xi^\pm) + 
\eps \cU'^{\pm} (\xi^\pm) + \ldots$, let us look for $\cU'^\pm$ as
\bea\label{E3}
\cU'^\pm(\xi^\pm) = c^{0\pm} \cY^{0}(\xi^\pm)
+c^{1\pm} \cY^{1 \pm}(\xi^\pm) + 
c^{2\pm} \cY^{2 \pm}(\xi^\pm).
\eea 
When $|\xi^{\pm}|\to+\infty$, we have
\begin{equation}\label{ExpansionInner}
u^\eps_-(x) = \left\{
\begin{array}{ll}
 w^+(y^\pm,-d)+\eps\,(c^{0\pm}+\dsp\frac{c^{1\pm}}{\pi} \ln \frac{1}{|\xi|} +c^{2\pm}\xi_1^{\pm} + O(|\xi|^{-1}))+\dots& \mbox{in }\bbR_+^2 \\[7pt]
v^0(\pm \ell )+\eps\,(c^{0\pm}+\dsp c^{1\pm} (\dsp\frac{\eta^{\pm}}{H(\pm\ell)}  + K^{1\pm } )+c^{2\pm}K^{2\pm } + O(e^{- \pi \eta^{\pm}/ H(\pm\ell)}))+\dots
&\mbox{in }\theta^{\pm}\setminus\bbR_+^2. 
\end{array}
\right.\hspace{-0.8cm}
\end{equation}
Comparing \ef{41}, \ef{44} and \ef{ExpansionInner}, we deduce that we ought to take 
\bea
c^{1 \pm} = \mp H( \pm \ell) \partial_s v^0(\pm \ell)\qquad \mbox{ and } \qquad c^{2 \pm } = \partial_y w^+(y^\pm , -d) . \label{E4}
\eea
Note that the constants $c^{0 \pm}$ in \ef{E3} still remain unfixed. Moreover, according to \ef{ExpansionInner}, the representation of $u^\eps_-(x)$ for $\xi^{\pm}\in\bbR_+^2$ such that $|\xi^{\pm}|\to+\infty$
includes the logarithmic term 
\bea
 \mp \eps \frac{H( \pm \ell)}{\pi} \ln \frac{1}{|\xi^\pm|}	 
\partial_sv^0(\pm \ell) . \label{E5}
\eea
Therefore, we will impose that the correction term in the outer expansion $u^\eps_- = w^+ + \eps u_-' + 
\ldots $ in $\Pi$ solves \ef{E7} and behaves as 
\bea\label{E8}
u'(x) = \mp  \frac{H( \pm \ell)}{\pi} \ln \frac{1}{r^\pm}	 
\partial_sv^0(\pm \ell) + \mrm{const.}
\eea
as $r^\pm \to 0^+$. One can prove the existence of a unique outgoing solution of this problem. We denote $\tilde{u}' (P^\pm)$ the corresponding constant in \ef{E8}. According to the definition of outgoing solutions (see \ef{DefOut}), we have the decomposition 
\bea\label{E9}
u' = \chi_-\,R'_-\,w^-+\chi_+\,T'\,w^+ + \ldots  
\eea
where the dots correspond to some term which is exponentially decaying at infinity in $\Pi$. Once $u'$ has been fixed, one obtains that the constants $c^{0 \pm}$ in \ef{E3} are given by
\bea\label{F1}
c^{0 \pm} =  \tilde{u}' (P^\pm) \pm  \frac{H( \pm \ell)}{\pi} 
\partial_sv^0(\pm \ell)  \ln \eps.  
\eea
Note that the term $\ln \eps$ appears because of the change of variable 
$ \xi^\pm \mapsto x$ in \ef{E5}. Continuing the matching procedure, one would write a problem similar to \ef{17}--\ef{18} for 
the correction term $\eps v'(s, \ln \eps)$ in the outer
expansion in the strip $\varpi^\eps$  and regard \ef{F1} as the data in this problem. This term depends linearly on $\ln \eps$, but it is not important in the sequel.

\BER\label{rem2.1}
Note that considering the stretching \ef{24} in the strip $\varpi^{\eps}$ leads for the variables $(s,n)$ to the transformation
\[
(s,n)\mapsto \eps^{-1}(s\mp l,n).
\]
Omitting the above assumption that $L$ is straight in a neighbourhood of $P^{\pm}$ would lead to a much more complicated change of variables 
\bea\label{mi3}
(s,n) \mapsto\eps^{-1} 
\Big(s\mp\ell, \frac{1}{H(s)} \big( H^+(\pm \ell)(n+\eps H^-(s))+H^-(\pm \ell)(n-\eps H^+(s))\big)\Big).
\eea
The Laplacian is still the main asymptotic part of the 
Helmholtz operator in the new variables \ef{mi3}, but the remainder is a second order differential operator with degenerating coefficients at the point $P^{\pm}$. This would make the estimation of the discrepancies in Section \ref{SectionJustif} much more cumbersome, cf. \cite{A4}. \ \ $\boxtimes$
\ENR

\subsection{Computing the scattering coefficients.} \label{sec3.3}
Once the singular solution $u_-'$ has been found, it is straightforward to
compute the coefficients $R_-'$, $T'$ in \ef{E9}. To this end, first we observe that the quantities 
\begin{equation}\label{EqnIndepD}
I^{\pm}=\sum_\pm \pm 
\int_{-d}^0 \Big( \overline{w^+(\pm D,z)} \partial_z u_-' (\pm D ,z) 
- u_-' (\pm D ,z) \overline{\partial_z w^+(\pm D,z)} \Big) dz
\end{equation}
are independent of $D>\ell_0$ (to show this property, integrate by parts and use that $u_-'$, $w^{\pm}$ satisfy the homogeneous Helmholtz equation in $\Pi$). Inserting the representation \ef{E9} in (\ref{EqnIndepD}) and taking the limit as $D\to+\infty$, we obtain 
\begin{equation}
I^-=2 i \sqrt{a^2 - \sfk^2 	} N\,R_-'\qquad\mbox{ and }\qquad I^+=2 i \sqrt{a^2 - \sfk^2 	} N\, T'.
\end{equation} 
Here $N$ is the same as in \ef{MainCorT}, such that 
\[
N = \Vert w^\pm (y, \cdot)  ; L^2(-d;0) \Vert^2 = (2a)^{-1}\,(e^{2ad}-e^{-2ad})+2d > 0.
\]
For $D>\ell_0$, $\delta>0$, define the domain  
\[
\Pi_\delta(D) := \Pi(D) \setminus  {\textstyle \bigcup\limits_{\pm}} 
{\textstyle\bigcup\limits_{j=1}^J}\,\overline{\bbB_\delta ( P_j^\pm)},
\]
with $\Pi(D) := \{ (y,z) \in \Pi :|y| < D \}$. Integrating by parts in  $\Pi_\delta(D)$ for a given $D>\ell_0$ and taking the limit $\delta\to0^+$, we also find 
\begin{equation}\label{F2}
\begin{array}{lcl}
I^{\pm}&=&\dsp\sum_\pm \sum_{j=1}^J 
\lim\limits_{\delta \to 0^+}  \delta
\int_0^\pi \Big( \overline{w^{\pm}(y,z)} \partial_{r_j^\pm} u_-' (y ,z) 
- u_-' (y ,z)   \overline{\partial_{r_j^\pm} w^{\pm}(y,z)}  \Big) \Big|_{r_j^\pm = \delta}d \varphi_j^\pm\\[8pt]
& =&- \dsp\sum_{j=1}^J   \Big( \overline{w^{\pm}(y_j^- ,-d)} H_j(-\ell_j) \partial_s v_j^0(-\ell_j) 
- \overline{ w^{\pm}(y_j^+, -d)}  H_j(\ell_j) \partial_s v_j^0(\ell_j) \Big).
\end{array}
\end{equation}
Here, $(r_j^\pm ,\varphi_j^\pm) \in \bbR_+ \times (0;\pi)$ is the polar 
coordinate system centered at $P_j^\pm$. Note that the logarithmic singularities in \ef{E8} have been taken into account to obtain the second line of \ef{F2}. Now using that $v^0_j(\pm \ell_j) = w^{+}(y^\pm_j,-d)$ (see \ef{18}), we get 
\begin{equation}\label{MainCorRBis}
\phantom{\mbox{and }}\hspace{0.61cm} R'_-=\cfrac{-i}{2N\sqrt{a^2-k^2}}\,\dsp\sum_{j=1}^J\int_{-\ell_j}^{\ell_j} H_j(s)\left(\bigg(\cfrac{dv^0_j}{ds}(s)\bigg)^2+\sfk^2(v^0_j(s))^2\right)\,ds
\end{equation}
\begin{equation}\label{MainCorTBis}
\mbox{and }\qquad T'=\cfrac{-i}{2N\sqrt{a^2-k^2}}\,\dsp\sum_{j=1}^J\int_{-\ell_j}^{\ell_j} H_j(s)\left(\bigg|\cfrac{dv^0_j}{ds}(s)\bigg|^2+\sfk^2|v^0_j(s)|^2\right)\,ds
\end{equation}
which are nothing but formulas \ef{MainCorR} and \ef{MainCorT}. 
As for the asymptotic expansions
\bea
R_-^\eps = \eps R_-' + \eps
\widetilde{R}_-^\eps\qquad\mbox{ and }\qquad T^\eps = 1 + \eps T' + \eps
\widetilde{T}^\eps,  \label{F4}
\eea
the following estimate will be proven in Section \ref{sec4.5}.

\BET\label{FF}
There are $\eps_0 > 0$ and $\delta \in (0;1/2)$ such that for $R'_-$, $T'$ respectively given by \ef{MainCorRBis}, \ef{MainCorTBis}, the remainders in \ef{F4} satisfy the estimate
\bea\label{F5}
| \widetilde{R}_-^\eps |+| \widetilde{T}_-^\eps | \leq c_\delta\,\eps^{1/2-\delta}.
\eea
Note that the constant $c_\delta$ depends on $\delta$ but not on $\eps \in (0;\eps_0]$. 
\ENT

\subsection{The global asymptotic approximation.} \label{sec3.4}

The matching method we have used up to now yields asymptotic expansions of $u^{\eps}_-$ in different zones. One may be interested in computing a global approximation of $u^{\eps}_-$ in the whole domain $\Pi^{\eps}$.  To obtain such an approximation, and more generally, to model boundary-value problems of type \ef{3}--\ef{5} on junctions of domains with different limit dimensions, an approach based on self-adjoint extensions of differential operators has been proposed in \cite{Naza96,na576}. Using this approach one can indeed model the initial problem by the system consisting of the equations
\ef{17}--\ef{18}, \ef{E7}--\ef{E9} and derive error estimates. Unfortunately we obtain limited accuracy and these estimate are not sufficient for the main goal of this paper.\\ 
\newline
Instead of using this model, we will work on the different expansions of $u^{\eps}_-$ and glue them to obtain a global approximation. The traditional approach \cite{Ilin} to do that, based on the use of partitions of unity, does not provide sufficient accuracy for the purpose of the paper. Instead, we will employ a trick with cut-off functions with overlapping supports, as introduced in \cite[Chap.\,2]{MaNaPl}. Note that the asymptotic structures we get with this technique have been shown to be equivalent with the method of compound asymptotic expansions (see \cite{ViLu,MaNaPl}).\\
\newline
The tunnel $\mathcal{T}^\eps_j$ enters the strip $\Pi$ through the two junction segments (see the bold line in Fig.\,\ref{fig1},\,b))
\bea\label{42}
\Xi^{\eps -}_j,\,\Xi^{\eps +}_j\quad\mbox{ such that }\quad\Xi^{\eps -}_j\cup\Xi^{\eps +}_j=\Gamma_d \cap \mathcal{T}^\eps_j\quad\mbox{ and }\quad P^{\pm}_j\in\Xi^{\eps \pm}_j.
\eea
Introduce smooth cut-off functions such that for $j=1,\dots,J$, we have 
\begin{equation}\label{CutOff1}
\cX_j^\pm=1\mbox{ in }\bbB_{R^\pm_j} (0) \qquad\mbox{ and }\qquad \cX_j^\pm=0\mbox{ outside }\bbB_{\widetilde{R}^\pm_j} (0).
\end{equation}
Then we define $X^{\eps}$ as 
\begin{equation}\label{CutOff2}
X^{\eps}(x)=1-\sum_{\pm}\sum_{j=1}^J\cX_j^{\eps\pm}(x)\qquad\mbox{ with }\qquad\cX_j^{\eps\pm}(x)=\cX_j^\pm(P^{\pm}_j+\eps^{-1}(x-P^{\pm}_j)).
\end{equation}
In \ef{CutOff1}, the radii $0<R^\pm_j<\widetilde{R}^\pm_j$ are chosen such that $\cX_j^{\eps \pm}=1$ in a neighbourhood of $\Xi^{\eps \pm}_j$. Finally we set $X^{\eps}_0:=X^{\eps}|_{\Pi}$ and $X^{\eps}_j:=X^{\eps}|_{\varpi^{\eps}_j}$, $j=1,\ldots,J$. Note that with a slight abuse of notation, we make no distinction between these functions and their extensions by zero to $\Pi^{\eps}$. Observe that $X^{\eps}_0$ is supported in $\overline{\Pi}$ while $X^{\eps}_j$ is supported in $\varpi^{\eps}_j$, $j=1,\ldots,J$.\\
\newline
Now we have everything to define our global approximation. In the straight part $\Pi$, set  
\begin{equation}\label{G1}
\begin{array}{lcl}
u_{\rm as}^\eps(x) & =& \dsp\phantom{+ }X_0^\eps(x) 
(w^+(x) + \eps u_-'(x))    \\
&&+\dsp\sum_\pm \sum_{j=1}^J \cX_j^\pm (x) ( \cU_j^{0\pm} + \eps 
\cU_j'^\pm (\xi_j^\pm))  \\
&&-\dsp\sum_\pm \sum_{j=1}^J X_0^\eps (x) \cX_j^\pm(x) 
\Big( \cU_j^{0\pm} +  \eps  \Big( c_j^{0 \pm} \ln \eps + c_j^{1 \pm}
\frac{1}{\pi} \ln \frac{1}{|\xi_j^\pm|} + c_j^{2 \pm}
\xi_{j1}^\pm \Big)\Big).
\end{array}
\end{equation}
In the curved thin strips $\varpi^{\eps}_j$, we set
\begin{equation}\label{G2}
\begin{array}{lcl}
u_{\rm as}^\eps(x) & =& \dsp\phantom{+ }X_j^\eps(x) v_j^0(s_j) \\

&&+\dsp\sum_\pm \sum_{j=1}^J \cX_j^\pm (x) \big( \cU_j^{0\pm} + \eps 
\cU_j'^\pm (\xi_j^\pm) \big) \\
 
&&-\dsp X_j^\eps (x) \sum_\pm \cX_j^\pm(x) 
\big( v_j(\pm \ell_j) + (s_j \mp \ell_j) \partial_s v_j^0( \pm \ell_j) \big). 
\end{array}
\end{equation}
These formulas need explanations. In $\Pi$,
the right hand side of \ef{G1} first contains the outer expansion \ef{MainAsymptoPi} multiplied by the cut-off function $X_0^\eps$ which annuls it in a neighbourhood of the junction segments \ef{42}. It also contains (second line) the inner expansions \ef{E1} multiplied by cut-off functions which annul them far from the junctions. Note that the terms matched in \S\ref{paragraphCorrection} are present in both expansions. However, their duplications are cancelled by the last subtrahend in \ef{G1} (which coincides exactly with these matched terms). A similar structure is found in \ef{G2}. However, the corrector $\eps v_j'( \ln \eps, s_j)$  of the outer expansion in $\varpi^{\eps}_j$ is not included to the first line on the right hand side of \ef{G2}. This fact simplifies the subtrahend in \ef{G2}, although it will lead to additional difficulties in the estimation of the discrepancies in Section \ref{sec4.5}. The authors do not know a simpler asymptotic structure which still provides asymptotically sharp error estimates.

\section{Justification of asymptotics.}\label{SectionJustif}
In this section, we prove error estimates to justify the asymptotic expansions derived formally in \S\ref{SectionAsymptotics}.
\subsection{Weighted Sobolev spaces with detached asymptotics.}\label{sec4.1}
For $\beta\in\bbR$, let $W_\beta^1 (\Pi^\eps)$ be the completion of the space 
$\mathscr{C}_c^\infty( \overline{\Pi^\eps} )$ with respect to the norm
\bea\label{52}
\Vert u^\eps ; W_\beta^1 (\Pi^\eps) \Vert := 
\big( \Vert e^{ \beta |y|} \nabla u^\eps; L^2(\Pi^\eps)\Vert^2 
+ \Vert e^{ \beta |y|}  u^\eps; L^2(\Pi^\eps)\Vert^2 
\big)^{1/2}. 
\eea
Note that this norm encloses an exponential weight. We start by considering the inhomogeneous water-wave problem 
\begin{equation}\label{54}
\begin{array}{|rcll}
- \Delta u^{\eps} + {\sf k}^2 u^{\eps} &=&  f^\eps &\mbox{ in }\Pi^\eps\\[3pt]
\partial_\nu u^{\eps} &=& 0 &\mbox{ on }\partial \Pi^\eps
\setminus \Gamma_0\\[3pt]
\partial_\nu u^{\eps} &=& \lambda u^{\eps} &\mbox{ on }\Gamma_0.\end{array}
\end{equation}
The weak formulation of \ef{54} writes
\bea\label{53}
\begin{array}{|l} 
\mbox{ Find }u^\eps\in W_{-\beta}^1 (\Pi^\eps)\mbox{ such that }\\[3pt]
\ (\nabla u^\eps , \nabla \psi )_{\Pi^\eps}
+ \sfk^2(u^\eps, \psi)_{\Pi^\eps} 
- \lambda (u^\eps ,\psi)_{\Gamma^0}
= F^\eps (\psi),\qquad \forall\psi\in W_\beta^1 (\Pi^\eps),
\end{array}  
\eea
where $F^\eps\in W_\beta^1(\Pi^\eps)^*$, the space of continuous antilinear functionals in
$W_\beta^1(\Pi^\eps)$. For $F^\eps$, for example, one can take  
\beas
F^\eps (\psi) = (f^\eps, \psi 
)_{\Pi^\eps}\qquad\mbox{ with }f^\eps\mbox{ such that } e^{-\beta|y|} f \in L^2(\Pi^\eps). 
\eeas
Using the Riesz representation theorem, define the continuous mapping $A_{-\beta}^\eps(\lambda):W_{-\beta}^1 (\Pi^\eps)\to W_\beta^1(\Pi^\eps)^*$ such that 
\begin{equation}\label{55}
\langle A_{-\beta}^\eps(\lambda)u^\eps, \psi\rangle_{\Pi^{\eps}}=( \nabla u^\eps , \nabla \psi )_{\Pi^\eps}
+ \sfk^2(u^\eps, \psi)_{\Pi^\eps} 
- \lambda (u^\eps ,\psi)_{\Gamma^0},\qquad \forall\psi\in W_\beta^1 (\Pi^\eps).  
\end{equation}
Here $\langle \cdot, \cdot\rangle_{\Pi^{\eps}}$ stands for the duality products between $W_\beta^1(\Pi^\eps)^*$ and $W_{\beta}^1 (\Pi^\eps)$. The Kondratiev theory \cite{Ko67} (see also, e.g., \cite[Chap.\,5]{NaPl}) yields the Fredholm property of the operators $A_{\pm \beta}^\eps(\lambda)$ if the following restriction holds for the weight exponent: 
\bea\label{57}
\beta \in (0;\sqrt{b^2+\sfk^2}). 
\eea
Here $b$ is the root of the transcendental equation \ef{TransEqb}. Following \cite{na642} (see also \cite[Chap.\,5]{NaPl}
and \cite{na489}), we introduce the weighted space with detached asymptotics
\bea
\bfW_\beta^1 (\Pi^\eps;\lambda ) \approx W_\beta^1(\Pi^\eps)
\oplus \bbC^2 \ni \{ \tilde{u}^\eps, a_\pm^\eps \}
\label{51}
\eea
consisting of functions of the form 
\bea\label{58}
u^\eps(x)  = \tilde{u}^\eps(x) + \sum_\pm\chi_\pm (y)
a_\pm^\eps w^\pm(x).
\eea
Here, $\chi_\pm$ are the cut-off functions introduced in \ef{CutOffChipm}. This space is endowed with the composite norm
\bea
\Vert u^\eps; \bfW_\beta^1(\Pi^\eps ; \lambda) \Vert :=
\Big( \Vert \tilde{u}^\eps; W_\beta^1(\Pi^\eps) \Vert^2+
\sum_\pm |a_\pm^\eps|^2 \Big)^{1/2} . \label{59}
\eea
Notice that the waves $w^\pm$ included into the decomposition 
\ef{58} depend of the spectral parameter $\lambda > \lambda_\dagger$, which is therefore included in the notation \ef{51}. We formulate an assertion proved in \cite[Thm.\,4.7]{na642}.

\BEL\label{lem4A}
Let $\lambda > \lambda_\dagger(\sfk)$ and assume that \ef{57} holds. Then the mapping
\bea
\bfA_\beta^\eps (\lambda) : \bfW_\beta^1 (\Pi^\eps ; \lambda)
\to W_{-\beta}^1(\Pi^\eps)^* \label{60}
\eea
is a Fredholm operator of index zero, and its null space satisfies the 
relation
\bea\label{61}
{\rm ker}\, \bfA_\beta^\eps (\lambda) = 
{\rm ker}\, A_\beta^\eps (\lambda).
\eea
\ENL 
\noindent A consequence of \ef{61} is that the kernel of $\bfA_\beta^\eps(\lambda)$ consists of trapped modes, that is solutions of the homogeneous problem \ef{54} which are exponentially decaying at infinity. From Lemma \ref{lem4A}, we also infer that the compatibility conditions 
\beas
F^\eps (\psi ) = 0,\qquad\forall\psi \in 
{\rm ker} \, A_\beta^\eps (\lambda) 
\eeas
imply the existence of a solution $u^\eps \in 
\bfW^1_\beta( \Pi^\eps; \lambda)$ to the problem \ef{53}. The orthogonality condition 
\bea
(u^\eps , \psi )_{\Pi^\eps} = 0,\qquad\forall\psi \in{\rm ker} \, A_\beta^\eps (\lambda) 
\eea
implies the uniqueness and the estimate 
\bea
\Vert u^\eps ; \bfW_\beta^1(\Pi^\eps ; \lambda) \Vert \leq c(\eps)
\Vert F^\eps; W_{-\beta}^1(\Pi^\eps)^* \Vert .\label{62}
\eea

\subsection{Limit problems in weighted spaces.}\label{sec4.2} 
We next present several known results concerning the limit problems. We denote by $W_{\beta,\gamma}^1 (\Pi)$ the Kondratiev space \cite{Ko67} which is obtained as the completion of $\mathscr{C}_c^\infty( \overline{\Pi} \setminus \{ P_1^\pm , \ldots,
P_J^\pm \})$ with respect to the norm
\bea \label{65}
\Vert \tilde{u} ; W_{\beta, \gamma}^1 (\Pi) \Vert := 
\big( \Vert e^{ \beta |y|} \rho_0^\gamma \nabla \tilde{u}; L^2(\Pi)\Vert^2 
+ \Vert e^{ \beta |y|} \rho_0^{\gamma-1} \tilde{u}; L^2(\Pi)\Vert^2 
\big)^{1/2}.
\eea
Here $\rho_0 (x) := \min \{ 1, r_1^\pm ,\ldots , r_J^\pm \}$ so that the
norm has two types of weights, the exponential one at infinity and 
the power weights near the points $P_1^\pm, \ldots, P_J^\pm$. The weights
make it possible to detach the asymptotics also near these points. Namely,
we define the space $\bfW_ {\beta, \gamma}^{1} (\Pi;\lambda)$ of functions of the form
\bea
u^0(x) = \tilde{u}^0(x) + \sum_\pm \Big( \chi_\pm (y)
a_\pm w^\pm(x) + \sum_{j=1}^J \cX_j^\pm (x) c_{j \pm}^0\Big)
\label{63}
\eea
where $\tilde{u}^0(x) \in W_{\beta, \gamma}^1(\Pi)$ (the cut-off functions $\cX_j^\pm$ are the ones introduced in \ef{CutOff1}). 
We choose the weight exponent
\bea
\gamma \in (-1;0) \label{64}
\eea
so that $\cX_j^\pm$ and $\cX_j^\pm \ln r_j^\pm$ do not belong
to $ W_{\beta, \gamma}^1(\Pi)$ while $\cX_j^\pm r_j^\pm$ does. Now, the Kondratiev theorem on asymptotics, see \cite{Ko67} and \cite[Chap.\,2]{NaPl},
leads to the following assertion. 

\BEL
\label{lem4B}
Let \ef{57}, \ef{64} be valid  and let $F^\eps \in W_{-\beta, -\gamma}^1
(\Pi)^* \subset W_{-\beta}^1(\Pi)^* $. Then, the unique solution $u^0
\in \bfW_\beta^1(\Pi;\lambda)$ of the problem 
\bea\label{uvF}
( \nabla u^0, \nabla \psi^0)_\Pi + \sfk^2 (u^0, \psi^0)_\Pi - \lambda(u^0, 
\psi^0)_{\Gamma_0} = F^0(\psi^0),\qquad \forall\psi^0 \in 
W_{-\beta}^1 (\Pi) , 
\eea
has the asymptotic form \ef{63} and satisfies the estimate
\bea
\Vert u^0 ; \bfW_{\beta, \gamma}^{1} ( \Pi; \lambda) \Vert &:=  & 
\Big( \Vert \tilde{u}^0; W_{\beta, \gamma}^1(\Pi) \Vert^2 +\sum_\pm \Big( |a_\pm^0|^2 + \sum_{j=1}^J |c_{j \pm}^0 | ^2 \Big)  \Big)^{1/2}
\leq c \Vert F^0 ; W_{-\beta, - \gamma}^1(\Pi)^* \Vert. \label{67}
\eea
\ENL
\noindent We emphasize that the unique solvability of the water-wave problem in the unperturbed strip $\Pi$, which is used in Lemma \ref{lem4B},
follows from a result similar to Lemma \ref{lem4A} for the problem set in $\Pi$ and from the fact that uniqueness can be obtained with the Fourier method.\\ 
\newline 
The next problem we consider is the mixed boundary value problem in the tunnel $\varpi_{j}^\eps$:
\bea
\big( \nabla V_j^\eps, \nabla \psi 
\big)_{\varpi_{j}^\eps } + \sfk^2 \big(  V_j^\eps, 
\psi  \big)_{\varpi_{j}^\eps}
= F_j^\eps(\psi),\qquad\forall \psi
\in H_0^1(\varpi_{j}^\eps; \Xi_j^{\eps\pm}) .
\label{68}
\eea
Here, $ H_0^1(\varpi_{j}^\eps; \Xi_j^{\eps\pm}) $
is the Sobolev space of functions satisfying the Dirichlet conditions on the junction segments $\Xi_j^{\eps\pm}$ defined in \ef{42}. Noticing that $d_t((h+t)^{-1})=-(h+t)^{-2}$ and using an integration by parts, one can establish the following classical one-dimensional Hardy inequality
\beas
\int_0^\ell (h+t)^{-2}  |\Psi(t)|^2 dt \leq 4\int_0^\ell
\Big| \frac{d \Psi}{dt } (t) \Big|^2 dt, \qquad\forall\Psi \in H^1(0;\ell) \ \mbox{ with } \  \Psi(0) =\Psi(\ell) = 0, 
\eeas
for all $ \ell > 0$, $h > 0$. This allows one to show that the standard Sobolev norm in 
$H_0^1(\varpi_{j}^\eps; \Xi_j^{\eps\pm})$ is equivalent
to the weighted norm
\bea
\ \ \ 
\Vert V^\eps ; H_0^1(\varpi_{j}^\eps; \Xi_j^{\eps\pm}) \Vert = 
\big( \Vert \nabla V^\eps ; L^2(\varpi_{j}^\eps) \Vert^2 
+  \Vert (\eps+ \rho_j)^{-1} V^\eps ; 
L^2(\varpi_{j}^\eps) \Vert^2 \big)^{1/2} \label{69}
\eea
where $\rho_j(s) = \min(|s-\ell_j|,|s+\ell_j|)$. This observation and the Riesz representation theorem imply the following simple assertion.

\BEL
\label{lem4C}
If $ F_j^\eps \in H_0^1(\varpi_{j}^\eps; 
\Xi_j^{\eps\pm})^*$, then problem \ef{68} has  a unique solution
$V_j^\eps  \in 
H_0^1(\varpi_{j}^\eps; \Xi_j^{\eps\pm})$
and there holds the estimate
\bea
\Vert V_j^\eps ; H_0^1(\varpi_{j}^\eps; 
\Xi_j^{\eps\pm})\Vert  
\leq 
C_j \Vert F_j^\eps; 
H_0^1(\varpi_{j}^\eps; \Xi_j^{\eps\pm})^* \Vert \label{70} 
\eea
where,  for some $\eps_j > 0$, the number $C_j$ is independent of 
$F_j^\eps$ and  $\eps \in (0;\eps_j]$.
\ENL

\noindent Finally, we consider the problem
\[
\begin{array}{|rlcl}
\  -\Delta  \cU &=& \cF& \quad\mbox{ in } \Theta^{\pm}_j \\[3pt]
\ \partial_{\nu} \cU &=& 0&\quad\mbox{ on } \partial\Theta^{\pm}_j 
\end{array}
\]
(see a picture of the unbounded domain 
$\Theta_j^\pm$ in Fig.\,\ref{fig3}) in the weighted space 
$\cW_{\tau,\sigma}^1 (\Theta_j^\pm)$ endowed with the norm
\bea\label{71}
\Vert \cU ; \cW_{\tau, \sigma}^1 (\Theta_j^\pm)\Vert  :=\big( \Vert  \cR_{\tau, \sigma} \nabla_\xi \cU ; L^2 (\Theta_j^\pm ) \Vert^2 
+\Vert  \cR_{\tau-1, \sigma}  \cU ; L^2 (\Theta_j^\pm ) \Vert^2 \big)^{1/2} ,
\eea
where 
\bea
\cR_{\tau, \sigma} (\xi) = (1 + |\xi|)^\tau \ \mbox{in} \ \bbR_+^2 , \quad\qquad
\cR_{\tau, \sigma} (\xi) = e^{\sigma \eta_j^{\pm}} \ \mbox{in} \
\theta_{j}^{\pm}\setminus\bbR_+^2. \label{72}
\eea

\BEL
\label{lem4D}
Let $\tau \in (-1;0)$, $\sigma \in (0;\pi/H_j(\pm \ell))$ and  
$\cF_j^\pm \in \cW_{-\tau,-\sigma}^1 (\Theta_j^{ \pm})^*$. Then, 
the problem 
\[
\big( \nabla_\xi \cU_j^\pm, \nabla_\xi \Psi_j^\pm 
\big)_{\Theta_j^\pm } = \cF_j^\pm(\Psi_j^\pm), \qquad\forall \Psi_j^\pm
\in \cW_{-\tau, -\sigma}^1 (\Theta_j^{\pm}) 
\]
has  a unique solution $\cU_j^\pm  \in 
\cW_{\tau, \sigma} (\Theta_j^\pm)$ and the estimate
\bea
\Vert \cU_j^\pm ; \cW_{\tau, \sigma}^1 (\Theta_j^\pm)\Vert  
\leq C_j^\pm  \Vert \cF_j^\pm; \cW_{-\tau, -\sigma}^1 (\Theta_j^\pm)^* \Vert \label{74} 
\eea
is valid with some coefficient $C_j^\pm> 0$  independent of $\cF_j^\pm$.
\ENL

\noindent This result again follows from the Kondratiev theory. According to a general result in \cite[Chap.\,5 and 6]{NaPl} (see also \S5 in the review paper \cite{na262}), the Neumann problem in the
domain $\Theta_j^\pm$ with two outlets at infinity always has a 
(certainly non-unique) solution $\cU_j^\pm \in \cW_{\tau, - \sigma}(
 \Theta_j^\pm)$ admitting the representation 
\beas
\cU_j^\pm (\xi) = \tilde{\cU}_j^\pm(\xi) + C_{j0}^\pm
+ C_{j1}^\pm H_j( \ell_j) ^{-1} \eta_j^\pm \ \ \ \mbox{in} \ \theta_{j}^\pm\setminus\bbR_+^2,
\eeas 
where $\tilde{\cU}_j^\pm$ decays exponentially as $\eta_j^\pm \to + \infty$.
Subtracting the linear combination $C_{j0}^\pm + C_{j1}^\pm \cY^{1\pm}(\xi)$  (see the definition of that functions in \ef{30})  makes the solution to belong to $\cW_{\tau , \sigma}^1
(\Theta_j^\pm)$, and it also becomes unique.

\subsection{Solvability of the problem in $\Pi^\eps$.}\label{sec4.3}
We define in $\Pi^\eps$ the weights 
\begin{equation}\label{75}
\begin{array}{ll}
 \rho_{\gamma,1}^\eps (x) = \left\{
\begin{array}{cl} 
(\eps  + \rho_0(x) )^\gamma &\mbox{ in }\Pi\\[3pt] 
\eps^\gamma &\mbox{ in }\varpi_{j}^\eps 
\end{array}\right.\qquad \mbox{ and }\qquad\rho_{\gamma,0}^\eps (x)=\left\{\begin{array}{cl} 
(\eps + \rho_0(x))^{\gamma-1}&\mbox{ in }\Pi\\[3pt] 
\eps^\gamma(\eps+\rho_j(x))^{-1}&\mbox{ in }\varpi_{j}^\eps
\end{array} \right.
\end{array} 
\end{equation}
by glueing the weights in $\Pi$ and in $\varpi_{j}^\eps$ of \S
\ref{sec4.2}. We denote by $W_{\beta,\gamma}^{1,\eps}(\Pi^\eps)$ the space $W_\beta^{1}(\Pi^\eps)$ endowed with the new norm
\bea\label{76}
\Vert u^\eps ; W_{\beta,\gamma}^{1,\eps}(\Pi^\eps)\Vert
:=\big( \Vert e^{\beta |y|} \rho_{\gamma, 1}^\eps \nabla
u^\eps ; L^2(\Pi^\eps)\Vert^2 + 
\Vert e^{\beta |y|} \rho_{\gamma, 0}^\eps 
u^\eps ; L^2(\Pi^\eps)\Vert^2 \big)^{1/2} . 
\eea
For any $\eps>0$, the norms of $W_{\beta,\gamma}^{1,\eps}(\Pi^\eps)$ and $W_\beta^1(\Pi^\eps)$ are equivalent but the constants of equivalence depend on $\eps$. Then we define the space with detached asymptotics $\bfW_ {\beta, \gamma}^{1, \eps} (\Pi^\eps)$ of functions of the form
\bea\label{77}
u^\eps(x) &=& \tilde{u}^\eps(x) + \sum_\pm \chi_\pm (y) 
a_\pm^\eps w^\pm(x) +\sum_{j=1}^J \Big( X_j^\eps(x) v_j^\eps(s_j) + \sum_\pm 
\cX_{j\pm}^\eps (x) c_{j \pm}^\eps \Big) ,  
\eea
where $\beta$ and $\gamma$ are fixed above,
\bea \label{78}
\tilde{u}^\eps \in W_{\beta, \gamma}^{1, \eps}
(\Pi^\eps),\quad a_\pm^\eps, c_{j \pm}^\eps \in \bbC, \quad v_j^\eps \in H^1 (- \ell_j;\ell_j) , \quad v_j^\eps(\pm \ell_j)
= c_{j \pm}^\eps, \ j= 1, \ldots, J.
\eea 
Note that in \ef{77}, the cut-off functions $X_j^\eps$, $\cX_{j\pm}^\eps$ have been defined in \ef{CutOff2}. Accordingly, the norm of this space is given by
\bea\label{79}
& & \Vert u^\eps  ; \bfW_{\beta, \gamma}^{1, \eps} 
(\Pi^\eps ; \lambda)   \Vert = \inf \Big( 
 \Vert \tilde{u} ^\eps  ; W_{\beta, \gamma}^{1, \eps} 
(\Pi^\eps)   \Vert + \sum_\pm |a_\pm^\eps|+\sum_{j=1}^J \Big( \sum_\pm |c_{j \pm}^\eps | +
\Vert v_j^\eps; H^1(-\ell_j;\ell_j) \Vert 
\Big)\Big) 
\eea
where the infimum is computed over all representations \ef{77} (observe that in \ef{77}, the representation of $u^\eps$ in $\Pi$ and in the tunnels $\varpi_{j}^\eps$ are not unique).

\BEL\label{lem4T}
The norms $\Vert\,\cdot\,; \bfW_\beta^1(\Pi^\eps ; \lambda) \Vert$ and $\Vert\,\cdot\,; \bfW_{\beta, \gamma}^{1, \eps} 
(\Pi^\eps ; \lambda)   \Vert$ are equivalent.
\ENL 
\noindent Proof. Let us prove the relations 
\bea \label{80}
c_\eps \Vert u^\eps  ; \bfW_{\beta}^{1} 
(\Pi^\eps ; \lambda)   \Vert
\leq 
\Vert u^\eps  ; \bfW_{\beta, \gamma}^{1, \eps} 
(\Pi^\eps ; \lambda)   \Vert
\leq C_\eps
\Vert u^\eps  ; \bfW_{\beta}^{1} 
(\Pi^\eps ; \lambda)   \Vert,  
\eea
where the constants $ c_\eps, C_\eps >  0$ should not depend on 
$u^\eps$ but surely depend on $\eps$. First, we verify the left inequality in \ef{80}. Let $u^\eps$ be an element of $\bfW_{\beta}^{1}(\Pi^\eps ; \lambda)$ with the decomposition
\begin{equation}\label{Repres1}
u^\eps(x)  = \tilde{u}^\eps(x) + \sum_\pm\chi_\pm (y)
a_\pm^\eps w^\pm(x).
\end{equation}
By definition of the $\inf$ in \ef{79}, there are $\tilde{v}^\eps \in W_{\beta, \gamma}^{1, \eps}
(\Pi^\eps)$, $c_{j \pm}^\eps \in \bbC$, $v_j^\eps \in H^1 (- \ell_j;\ell_j)$ with $v_j^\eps(\pm \ell_j)= c_{j \pm}^\eps$, $j= 1, \ldots, J$ such that 
\begin{equation}\label{Repres2}
u^\eps(x) = \tilde{v}^\eps(x) + \sum_\pm \chi_\pm (y) 
a_\pm^\eps w^\pm(x) +\sum_{j=1}^J \Big( X_j^\eps(x) v_j^\eps(s_j) + \sum_\pm \cX_{j\pm}^\eps (x) c_{j \pm}^\eps \Big)
\end{equation}
and such that 
\begin{equation}\label{defInf}
\Vert \tilde{v} ^\eps  ; W_{\beta, \gamma}^{1, \eps} 
(\Pi^\eps)   \Vert + \sum_\pm |a_\pm^\eps|+\sum_{j=1}^J \Big( \sum_\pm |c_{j \pm}^\eps | +
\Vert v_j^\eps; H^1(-\ell_j;\ell_j) \Vert\Big) \le 2\,\Vert u^\eps  ; \bfW_{\beta, \gamma}^{1, \eps}(\Pi^\eps ; \lambda)\Vert.
\end{equation}
Then we can write
\begin{equation}\label{estimEq1}
\begin{array}{lcl}
\hspace{-0.3cm}\dsp\Vert u^\eps  ; \bfW_{\beta}^{1} 
(\Pi^\eps ; \lambda)   \Vert =\dsp \Vert \tilde{u} ^\eps  ; W_{\beta}^{1} 
(\Pi^\eps)   \Vert + \sum_\pm |a_\pm^\eps| & \le &  \dsp\Vert \tilde{v} ^\eps  ; W_{\beta}^{1} 
(\Pi^\eps)   \Vert + \sum_\pm |a_\pm^\eps|+\dsp\Vert \tilde{v} ^\eps-\tilde{u} ^\eps  ; W_{\beta}^{1} (\Pi^\eps)\Vert\\
& \le &  c_{\eps}\dsp\Vert \tilde{v} ^\eps  ; W_{\beta, \gamma}^{1, \eps} 
(\Pi^\eps)   \Vert + \sum_\pm |a_\pm^\eps|+\dsp\Vert \tilde{v} ^\eps-\tilde{u} ^\eps  ; W_{\beta}^{1} (\Pi^\eps)\Vert\\
& \le &  2\,c_{\eps}\,\Vert u^\eps  ; \bfW_{\beta, \gamma}^{1, \eps}(\Pi^\eps ; \lambda)\Vert+\dsp\Vert \tilde{v} ^\eps-\tilde{u} ^\eps  ; W_{\beta}^{1} (\Pi^\eps)\Vert.
\end{array}
\end{equation}
Now comparing the two representations \ef{Repres1}, \ef{Repres2} and using again \ef{defInf}, one finds 
\begin{equation}\label{estimEq2}
\dsp\Vert \tilde{v} ^\eps-\tilde{u} ^\eps  ; W_{\beta}^{1} (\Pi^\eps)\Vert \le c_{\eps}\sum_{j=1}^J \Big( \sum_\pm |c_{j \pm}^\eps | +
\Vert v_j^\eps; H^1(-\ell_j;\ell_j) \Vert\Big) \le 2\,c_{\eps}\Vert u^\eps  ; \bfW_{\beta, \gamma}^{1, \eps}(\Pi^\eps ; \lambda)\Vert.
\end{equation}
Inserting \ef{estimEq2} in \ef{estimEq1} leads to the left inequality of \ef{80}. \\
\newline
The right inequality of \ef{80} is proven by setting
$\tilde{u}^\eps + \sum_\pm \chi_\pm a_\pm^\eps w^\pm = 
u^\eps$ and $v_j^\eps= 0$, 
$c_{j\pm}^\eps = 0$, and referring again to the obvious equivalence of the norms of $W_{\beta,\gamma}^{1,\eps}(\Pi^\eps)$ and $W_\beta^1(\Pi^\eps)$. We emphasize again that the constants related to this 
equivalence relation depend on $\eps$. That explains the very difference in norming the same function space.\ \ $\boxtimes$\\
\newline
The weak formulation of problem \ef{54} in $\bfW_{\beta, \gamma}^{1,\eps}  (\Pi^\eps ; \lambda) $ is similar to \ef{53}. It allows one to define the bounded linear operator
\bea
\bfA_{\beta, \gamma}^{\eps}  ( \lambda) : 
\bfW_{\beta, \gamma}^{1,\eps}  (\Pi^\eps ; \lambda) 
\to
W_{-\beta, -\gamma}^{1,\eps}  (\Pi^\eps )^* ,  \label{81} 
\eea
which coincides with \ef{60}, although it is related to different norms. In the next section, we construct an approximate inverse  operator
\bea
\bfR_{\beta, \gamma}^{\eps}  ( \lambda) : 
 W_{-\beta,- \gamma}^{1,\eps}  (\Pi^\eps)^*  \to 
\bfW_{\beta, \gamma}^{1,\eps}  (\Pi^\eps ; \lambda)  \label{82}
\eea
such that the difference
\bea
{\rm Id} - \bfA_{\beta, \gamma}^{\eps}  ( \lambda) 
\bfR_{\beta, \gamma}^{\eps}  ( \lambda) :
W_{-\beta,- \gamma}^{1,\eps}  (\Pi^\eps)^* 
\to  W_{-\beta,- \gamma}^{1,\eps}  (\Pi^\eps)^*   \label{83}
\eea
has a small operator norm $o(1)$ as $\eps \to 0^+$. In view of the classical result concerning Neumann
series, this shows that $\bfA_{\beta, \gamma}^{\eps}  ( \lambda) \bfR_{\beta, \gamma}^{\eps}  ( \lambda)$ is invertible for $\eps$ small enough. We deduce that $\bfA_{\beta, \gamma}^{\eps}(\lambda)^{-1}:=\bfR_{\beta, \gamma}^{\eps}  ( \lambda)
\big( \bfA_{\beta, \gamma}^{\eps}  ( \lambda) 
\bfR_{\beta, \gamma}^{\eps}  ( \lambda) \big)^{-1}$ is  a  continuous right inverse of $\bfA_{\beta, \gamma}^{\eps}(\lambda)$. From the fact that $\bfA_{\beta, \gamma}^{\eps}(\lambda)$ is Fredholm of index zero (Lemma \ref{lem4A}), we infer that $\bfA_{\beta, \gamma}^{\eps}(\lambda)$ is actually an isomorphism and that $\bfA_{\beta, \gamma}^{\eps}(\lambda)^{-1}$ is its true inverse for $\eps$ small enough. Thus, we can state the following important result which yields in particular a stability estimate. 

\BET
\label{th4E}
Set $\beta\in(0;\sqrt{b^2+\sfk^2})$ and $\gamma\in(-1/2;0)$ where $b$ is defined in (\ref{TransEqb}). There is some $\eps_ 0>0$ such that for all $\eps \in (0;\eps_0]$, for all $F^\eps \in W_{-\beta, - \gamma}^{ 1 , \eps}
(\Pi^\eps)^* $, the problem \ef{54} in $\Pi^{\eps}$ has a unique solution
$u^\eps \in\bfW_{\beta, \gamma}^{1,\eps}  (\Pi^\eps ; \lambda)$. In particular, trapped modes are absent. Moreover,  the estimate
\bea
\Vert u^\eps ; \bfW_{\beta, \gamma}^{1,\eps}  (\Pi^\eps ; 
\lambda) \Vert  \leq C \Vert F^\eps; 
W_{-\beta, -\gamma}^{1,\eps}  (\Pi^\eps ; \lambda)^* \Vert
\label{84}
\eea
is valid with a constant $C>0$ independent of $F^\eps$ and, what is of the most importance, of $\eps \in (0;\eps_0]$.
\ENT

\subsection{Construction of an approximate inverse operator.} \label{sec4.4}
We will find the operator \ef{82} as the sum
\bea
\bfR_{\beta, \gamma}^{\eps }  ( \lambda) = 
\bfR_{\beta, \gamma}^{\eps 0}  ( \lambda)
+ \sum_{j=1}^J \Big( \bfR_{\beta, \gamma}^{\eps j} 
+ \sum_\pm \bfR_{\beta, \gamma}^{\eps j \pm} \Big), \label{85}
\eea
where the terms will be defined below. We fix a functional
$F^\eps \in W_{-\beta, -\gamma}^{1, \eps } (\Pi^\eps)^* $
and denote its norm by $N^\eps$.\\
\newline
{\bf Step 1.} Denote $u^{\eps0} \in
\bfW_{\beta, \gamma}^{1}  (\Pi;\lambda)$ the solution of problem \ef{uvF}, which is well-defined according to Lemma \ref{lem4B}, with the
right-hand side $F^{\eps 0}$ such that 
\beas
F^{\eps 0} (v^0) = F^\eps(X_0^\eps\,v^0),\qquad \forall v^0\in W_{-\beta,-\gamma}^1(\Pi).
\eeas 
Here $X_0^\eps$ is the cut-off function introduced after \ef{CutOff2} which vanishes at the points $P_j^\pm$. According to the definitions \ef{76} and \ef{75}, we have $\|X_0^\eps\,v^0;W_{-\beta, -\gamma}^{1, \eps } (\Pi^\eps)\|\le c\,\|v^0;W_{-\beta, -\gamma}^{1} (\Pi)\|$. We deduce the estimate
\bea 
\Vert F^{\eps 0 } ; W_{-\beta,-\gamma}^1(\Pi)^* \Vert \leq cN^\eps.
\label{86}
\eea
Note that in \ef{86} the constant $c>0$ is independent of $\eps$. We emphasize that in what follows, if a constant depends on $\eps$, then we will write it explicitly. Thus, the estimate \ef{67} with the bound \ef{86} holds for the ingredients $\tilde{u} ^{\eps 0}$,
$a_\pm^{\eps 0}$ and $c_{j \pm}^{\eps 0}$ 
of the representation \ef{63}.  We set
\begin{equation}\label{87}
\bfR_{\beta, \gamma}^{\eps 0}  ( \lambda) F^\eps :=  \bfu_0^{\eps} :=\dsp X_0^\eps(x) \tilde{u}^{\eps0} (x)  
+ \sum_\pm \chi_\pm(y) a_\pm^{\eps 0} w^\pm (x) +\sum_{j=1}^J \Big( X_j^\eps(x) v_j^{\eps0} (s_j) 
+ \sum_\pm \cX_{j\pm}^\eps (x) c_{j \pm}^{\eps0} \Big),
\end{equation}
where $v_j^{\eps 0}$ is the unique solution of the equation \ef{17}
with the Dirichlet data $v_j^{\eps 0} (\pm \ell_j)
= c_{j \pm}^{\eps 0}$. Note in particular that 
\[
\Vert v_j^{\eps 0} ; H^2 (-\ell_j;\ell_j) \Vert \leq c 
\sum_\pm |c_{j \pm}^{\eps 0}|
\leq c N^\eps. 
\]
From the definition of the norm on $\bfW_{\beta, \gamma}^{1 , \eps }  ( \Pi^\eps;\lambda)$ in \ef{79}, we deduce that  
\bea
\Vert \bfR_{\beta, \gamma}^{\eps 0}  ( \lambda) F^\eps ;
\bfW_{\beta, \gamma}^{1 , \eps }  ( \Pi^\eps;\lambda)\Vert
\leq c_0N^\eps . \label{89}
\eea
\noindent {\bf Step 2.} Let $V_j^\eps$ be the solution of the mixed boundary value problem \ef{68} with the right-hand side $F_j^\eps$ such that 
\bea\label{90}
F_j^\eps(\psi) = F^\eps(X_j^\eps
\,\psi),\qquad\forall  \psi\in H_0^1(\varpi_{j}^\eps; \Xi_j^{\eps\pm}) .
\eea
Note that $V_j^\eps$ is well-defined according to Lemma \ref{lem4C} and we have 
\beas
\Vert V_j^\eps ; H_0^1(\varpi_{j}^\eps; 
\Xi_j^{\eps\pm})\Vert  
\leq 
c \Vert F_j^\eps; 
H_0^1(\varpi_{j}^\eps; \Xi_j^{\eps\pm})^* \Vert .
\eeas
Moreover, using \ef{70}
and \ef{75}, \ef{69}, we can write 
\[
\begin{array}{lcl}
 |F_j^\eps(\psi)| &\leq& 
\Vert F^\eps; W_{-\beta, -\gamma}^{1, \eps} (\Pi^\eps)^*  
\Vert \, 
\Vert X_j^\eps \psi; W_{-\beta, -\gamma}^{1, \eps} (\Pi^\eps) \Vert
\\[6pt]
& \le & cN^\eps \eps^{-\gamma} 
\big( \Vert \nabla (  X_j^\eps \psi) ;
L^2(\varpi_{j}^\eps) \Vert + 
\dsp\Vert (\eps  + \rho_j)^{-1} X_j^\eps 
\psi ; L^2(\varpi_{j}^\eps) \Vert \big) 
\\[4pt]
&\le & cN^\eps \eps^{-\gamma} 
\big( \Vert \nabla  \psi ;  L^2(\varpi_{j}^\eps) \Vert + 
\Vert \psi \nabla   X_j^\eps  ;
L^2(\varpi_{j}^\eps) \Vert + 
\dsp\Vert (\eps  + \rho_j)^{-1} 
\psi ; L^2(\varpi_{j}^\eps) \Vert \big) \\[4pt]
&\le & cN^\eps \eps^{-\gamma} \Vert \psi ;
H_0^1 (\varpi_{j}^\eps ; \Xi_j^{\eps \pm} ) \Vert.  
\end{array}
\]
We deduce that $\Vert F_j^\eps ; H_0^1 (\varpi_{j}^\eps ; \Xi_j^{\eps \pm} )^{\ast} \Vert \leq c N^\eps \eps^{-\gamma}$. Here, we have taken into account that $\nabla X_j^\eps =0$
for $s_j \in [- \ell_j+ \widetilde{R}_j^- \eps ;\ell_j -\widetilde{R}_j^+\eps ]$ so that $|\nabla X_j^\eps
(x) | \leq c\big( \eps  + \rho_j(x) \big)^{-1}$ in $\varpi_{j}^\eps$. Therefore, setting 
\bea\label{91}
\bfR_{\beta,\gamma}^{\eps j} F^\eps:=\bfu_j^\eps:=X_j^\eps(x) V_j^\eps(x), 
\eea
we find 
\bea\label{91bis}
\Vert  \bfR_{\beta,\gamma}^{\eps j} F^\eps ; 
\bfW_{\beta,\gamma}^{1,\eps} (\Pi^\eps; \lambda)\Vert
\leq c N^\eps
\eea
(take a representation with $v_j^{\eps}=0$ in (\ref{79}) and use the definition of the norm of $W_{\beta, \gamma}^{1, \eps}
(\Pi^\eps)$.\\[4pt]
{\bf Step 3.} Let us calculate the discrepancy (mismatch with respect to the data $F^\eps$) left in the problem \ef{53} by the sum
\begin{equation}\label{DecompoSum}
\bfu^\eps = \bfu_0^\eps + \bfu_1^\eps +
\ldots + \bfu_J^\eps .
\end{equation}
We denote $\bfu_0^{\eps j} = \bfu_0^\eps \big|_{\varpi_{j}^\eps}$. For every $\psi \in W_{- \beta, -\gamma}^1(\Pi^\eps)$ with compact support (we can work by density), we find
\[
\big( \nabla \bfu_0^{\eps},\nabla \psi \big)_{\Pi^{\eps}}
+ \sfk^2 \big(  \bfu_0^{\eps} ,  \psi \big)_{\Pi^{\eps}}
- \lambda \big(\bfu_0^{\eps} ,  \psi \big)_{\Gamma_0}=
F^\eps (X_0^\eps\psi) + \bff_0^\eps(\psi)+\sum_{j=1}^J\bff_0^{\eps j}(\psi)
\]
\begin{equation}\label{92}
\mbox{ with }\qquad 
\begin{array}{|lcl}
\,\bff_0^\eps(\psi) &=& - \big( \nabla \big( (1- 
X_0^\eps) \tilde{u}^{\eps0 }\big) , \nabla 
\psi \big)_\Pi
- \sfk^2 \big( ( 1-  X_0^\eps) \tilde{u}^{\eps0 } , 
\psi \big)_\Pi \\[7pt]
\,\bff_0^{\eps j} (\psi) &=& \big( \nabla \bfu_0^{\eps j}
, \nabla\psi\big)_{\varpi_{j}^\eps}
+ \sfk^2 \big( \bfu_0^{\eps j} ,\psi \big)_{
\varpi_{j}^\eps}.
\end{array}
\end{equation}
In the following, we shall say that a functional $f^\eps$ is $O(\eps^\alpha)$ if it satisfies $\Vert f^\eps ; W_{-\beta,-\gamma}^{1, \eps} 
(\Pi^\eps)^\ast \Vert\leq c N^\eps \eps^\alpha$. Note that $\bff_0^\eps$  is a sum of terms located in $\eps$-neighbourhoods of the
points $P_1^\pm, \ldots, P_J^\pm$. One can prove that it is $O(\eps^0)$. We will compensate it by boundary layers in the next step. In contrast, we prove now that the $\bff_0^{\eps j}$ are small as $\eps$ tends to zero.\\
\newline
Observe that the $L^2$-norm of $\psi\in W_{-\beta,-\gamma}^{1, \eps}(\Pi^\eps)$ has the weights of order 
$\eps^{-\gamma-1}$ in the vicinity of the points $P_j^\pm$
and of order $\eps^{-\gamma}$ in $\varpi_{j}^\eps$, see
\ef{75} and \ef{76}. Note also that $\bfu_0^{\eps j}$ is constant in $\bbB_{R^\pm_j} (P_j^\pm)$. Recalling the formulas \ef{12}--\ef{14} and the relation $|A(n,s) -1|\leq c\,\eps $ in $\varpi_{j}^\eps$, we find
\begin{equation}\label{93}
\begin{array}{l}
\bff_0^{\eps j} (\psi) = \dsp\big( \nabla \bfu_0^{\eps j},\nabla \psi\big)_{\varpi_{j}^\eps}+ \sfk^2 \big( \bfu_0^{\eps j} ,\psi \big)_{
\varpi_{j}^\eps} = \big( - \Delta  \bfu_0^{\eps j} + \sfk^2  \bfu_0^{\eps j}
, \psi\big)_{\varpi_{j}^\eps}
+ \big( \partial_\nu  \bfu_0^{\eps j} ,  
\psi\big)_{\partial \varpi_{j}^\eps} \\[10pt]
\hspace{0.7cm}=\dsp\big(- \Delta  \bfu_0^{\eps j} + \sfk^2  \bfu_0^{\eps j}
, \psi\big)_{\varpi_{j}^\eps}-\sum_\pm \eps \int_{-\ell_j + \eps R^{-}_j}^{
\ell_j - \eps R_j^{+}}  \frac{dH_j^\pm}{ds} (s) 
\frac{d \bfu_0^{\eps j}}{ds} (s)\overline{\psi}\big( \pm \eps
H_j^\pm (s) ,s \big )\,ds + O(\eps^{\gamma +3/2}).
\end{array}
\end{equation}
Now using \ef{87}, in $\varpi_{j}^\eps$ we obtain
\[
- \Delta  \bfu_0^{\eps j} + \sfk^2  \bfu_0^{\eps j}=\sum_{\pm}2\nabla
v_j^{\eps0} (s_j)\cdot\nabla \cX_{j\pm}^\eps+(v_j^{\eps0} (s_j)-c_{j \pm}^{\eps0})\Delta\cX_{j\pm}^\eps+\sfk^2c_{j \pm}^{\eps0}\cX_{j\pm}^\eps+ O(\eps^{\gamma +3/2}).
\]
Using that $v_j^{\eps0}\in H^2(-\ell_j;\ell_j)\subset\mathscr{C}^1[-\ell_j;\ell_j]$ with the estimate
\[
\sup_{[-\ell_j;\ell_j]}|v_j^{\eps0}(s)|+\bigg|\frac{dv_j^{\eps0}}{ds}(s)\bigg| \le c\Vert v_j^{\eps 0} ; H^2 (-\ell_j;\ell_j) \Vert \leq c 
\sum_\pm |c_{j \pm}^{\eps 0}|
\leq c N^\eps,
\]
one can write
\begin{equation}\label{EstimInter1}
|\big(- \Delta  \bfu_0^{\eps j} + \sfk^2  \bfu_0^{\eps j}
, \psi\big)_{\varpi_{j}^\eps}| \le  c N^\eps\eps^{\gamma+1}.
\end{equation}
On the other hand, one can show the estimate 
\begin{equation}\label{EstimInter2}
\Big|\sum_\pm \eps \int_{-\ell_j + \eps R^{-}_j}^{
\ell_j - \eps R_j^{+}}  \frac{dH_j^\pm}{ds} (s) 
\frac{d \bfu_0^{\eps j}}{ds} (s)\overline{\psi}\big( \pm \eps
H_j^\pm (s) ,s \big )\,ds \Big| \le c N^\eps\eps^{\gamma+1/2}.
\end{equation}
Inserting \ef{EstimInter1} and \ef{EstimInter2} in \ef{93}, finally we obtain 
\begin{equation}\label{94}
|\bff_0^{\eps j} (\psi)| \leq cN^\eps 
\eps^{\gamma+1/2} \Vert \psi ; 
W_{- \beta, -\gamma}^{1, \eps} (\Pi^\eps) \Vert. 
\end{equation}
Now let us turn to the discrepancy left by the terms $\bfu_j^\eps$ in the  sum \ef{DecompoSum}. These functions vanish at the segments $\Xi_j^{\eps \pm }$ and we have
\[
\big( \nabla \bfu_j^\eps, \nabla \psi\big)_{
\varpi_{j}^\eps} + \sfk^2 (  \bfu_j^\eps,  
\psi )_{\varpi_{j}^\eps}
= F^\eps (X_j^\eps \psi) 
+ \bff_j^\eps(\psi)
\]
\bea
\mbox{with }\qquad
\bff_j^\eps (\psi) = 
- \big( \nabla((1- X_j^\eps)V_j^\eps), 
\nabla \psi \big)_{\varpi_{j}^\eps}
-\sfk^2  \big( (1- X_j^\eps)V_j^\eps), 
 \psi \big)_{\varpi_{j}^\eps} .\label{95}
\eea
Finally, since $F^\eps$ is antilinear, we get 
\begin{equation}\label{decompoRemainder}
\begin{array}{ll}
&\big(\nabla \bfu^{\eps},\nabla \psi \big)_{\Pi^{\eps}}
+ \sfk^2 \big(  \bfu^{\eps} ,  \psi \big)_{\Pi^{\eps}}
- \lambda \big(\bfu^{\eps} ,  \psi \big)_{\Gamma_0}\\[2pt]
=&\dsp F^\eps (\psi)+\Big(\bff_0^\eps(\psi)+\sum_{j=1}^J\big(\bff_j^\eps(\psi)-\sum_{\pm}F^\eps(\cX_j^{\eps\pm}\psi)\big)\Big) +\dsp\sum_{j=1}^J\bff_0^{\eps j}(\psi).
\end{array}
\end{equation}
Note that to obtain the identity, we used in particular the sequence of equalities 
\[
\sum_{j=0}^JF^\eps (X_j^\eps\psi)=F^\eps (X^\eps\psi)=F^\eps (\psi)-\sum_{\pm}\sum_{j=1}^JF^\eps(\cX_j^{\eps\pm}\psi).
\]
The second element of the right hand side of \ef{decompoRemainder} is supported in $\eps$-neighbourhoods of the points $P_1^\pm, \ldots, P_J^\pm$. Now we compensate it by boundary layers terms.\\
\newline 
{\bf Step 4.} We take an arbitrary $\Psi_{j}^{\pm} \in \cW_{-\tau, -\sigma}^1(
\Theta_j^\pm)$ and in order to compensate the second term of the right hand side of \ef{decompoRemainder}, we introduce the functionals $\cF_{j}^{\eps\pm}$ such that
\begin{equation}\label{97}
\cF_{j}^{\eps\pm} (\Psi_{j}^{\pm} ) = -\Big(\bff_0^\eps(\psi_{j}^{\pm})+\bff_j^\eps(\psi_{j}^{\pm})-F^\eps(\cX_j^{\pm}\psi_{j}^{\pm})\Big),
\end{equation}
with  $\psi_{j}^{\pm} (x) = \Psi_j^{\eps\pm} \big(P_j^\pm + \eps\xi_j^\pm \big)$. We remind the reader that $\xi_j^\pm$ are the stretched coordinates defined in \ef{24}. Due to the presence of the compactly supported cut-off functions in \ef{97}, the map $\cF_{j}^{\eps\pm}$ is continuous in $\cW_{-\tau, -\sigma}^1( \Theta_j^\pm)$ for any weight indices $\tau$ and 
$\sigma$, in particular for
\bea\label{98}
\tau \in [\gamma,0)\quad\qquad \mbox{and} \quad\qquad\sigma \in (0;\pi /H_j(\ell_j)).
\eea
Again, observe that the  weights for $\nabla_x \psi_{j}^{\pm}$ and 
$\psi_{j}^{\pm}$ in the norm of $W_{-\beta, - \gamma}(\Pi^\eps)$, 
are of the orders $\eps^{-\gamma}$ and $\eps^{-\gamma-1}$,
respectively, in the vicinity of the point $P_j^\pm$, see \ef{75} and
\ef{76}. Hence, we  conclude that
\begin{equation}\label{EstimEquiv}
c \big\Vert \Psi_{j}^{\pm} ; \cW_{-\tau, -\sigma}^1\big(
\Theta_j^\pm \cap \bbB_R(0)\big)\big\Vert 
\leq \eps^\gamma 
\big\Vert \psi_{j}^{\pm} ; W_{-\beta, -\gamma}^1\big(
\Pi^\eps \cap \bbB_{R \eps}(P_j^\pm \big)\big\Vert \le C\big\Vert \Psi_{j}^{\pm} ; \cW_{-\tau, -\sigma}^1\big(
\Theta_j^\pm \cap \bbB_R(0)\big)\big\Vert ,
\end{equation}
where the constants $C \geq c > 0$ depend on $R>0$ but not on $\Psi_{j}^{\pm}$, $\eps \in (0;\eps_0]$. Using the definitions of the three terms on the right-hand side of \ef{97}, we get 
\bea\label{Fe}
\begin{array}{lcl}
\big| \cF_{j}^{\eps\pm} (\Psi_{j}^{\pm})\big|
\leq 
cN^\eps 
\big\Vert \psi_{j}^{\pm} ; W_{-\beta, -\gamma}^1\big(
\Pi^\eps \cap \bbB_{R \eps}(P_j^\pm ) \big)\big\Vert &\le &
cN^\eps \eps^{-\gamma} \big\Vert \Psi_{j}^{\pm} ; \cW_{-\tau, -\sigma}^1(\Theta_j^\pm\cap \bbB_R(0)) \Vert\\[4pt]
&\le &
cN^\eps \eps^{-\gamma} \big\Vert \Psi_{j}^{\pm} ; \cW_{-\tau, -\sigma}^1(\Theta_j^\pm)\Vert.
\end{array}
\eea
As a consequence, for any $\tau$, $\sigma$ satisfying (\ref{98}), Lemma \ref{lem4D} guarantees the existence of a unique function $\cU_{j}^{\eps\pm}\in \cW_{\tau, \sigma}^1(\Theta_j^\pm)$ satisfying
\[
\big( \nabla_\xi \cU_{j}^{\eps\pm}, \nabla_\xi \Psi_j^\pm 
\big)_{\Theta_j^\pm } = \cF_j^\pm(\Psi_j^\pm), \qquad\forall \Psi_j^\pm
\in \cW_{-\tau, -\sigma}^1 (\Theta_j^{\pm}). 
\]
Moreover, we have the estimate
\bea\label{FeN}
\Vert \cU_{j}^{\eps\pm} ; \cW_{\tau, \sigma}^1(\Theta_j^\pm) \Vert
\leq C_{\tau, \sigma}\,\eps^{-\gamma} N^\eps. 
\eea
We then set
\bea\label{99}
\bfR_{\beta, \gamma}^{\eps j \pm} F^\eps (x) =\mathfrak{X}_j^\pm(x)\mathfrak{u}_{j}^{\eps\pm}(x)\qquad\mbox{ with }\qquad 
\begin{array}{|ll} 
 \mathfrak{X}_j^\pm(x)=\cX_j^\pm(x- P_j^\pm)\\[3pt]
 \mathfrak{u}_{j}^{\eps\pm}(x)=\cU_{j}^{\eps\pm}\big(\eps^{-1} (x- P_j^\pm)\big).
\end{array}
\eea
Let us fix $\tau$, $\sigma$ such that \ef{98} holds. Working as in (\ref{EstimEquiv}) and taking into account \ef{Fe}, we can write
\begin{equation}\label{X1}
\Vert \bfR_{\beta, \gamma}^{\eps j \pm} F^\eps; \bfW_{\beta, \gamma}^{1, \eps} (\Pi^\eps; \lambda)\Vert 
\leq c \Vert  \mathfrak{X}_j^\pm \mathfrak{u}_{j}^{\eps\pm}
; W_{\beta, \gamma}^{1,\eps} (\Pi^\eps)\Vert 
\le c \eps^\gamma  \Vert  \cU_{j}^{\eps\pm}; \cW_{\gamma, \sigma}^1(\Theta_j^\pm)\Vert 
 \leq c_{\gamma, \sigma} N^\eps .
\end{equation}
{\bf Conclusion.} Now we have completed the construction of the operator $\bfR_{\beta, \gamma}^{\eps }:W_{-\beta,- \gamma}^{1,\eps}  (\Pi^\eps)^*  \to 
\bfW_{\beta, \gamma}^{1,\eps}  (\Pi^\eps ; \lambda)$ in \ef{85}. It is the sum of the operators $\bfR_{\beta, \gamma}^{\eps 0}$, $\bfR_{\beta,\gamma}^{\eps j}$, $\bfR_{\beta, \gamma}^{\eps j \pm}$ defined respectively in \ef{87}, \ef{91}, \ef{99}. Its operator norm is bounded uniformly with respect to $\eps \in (0;\eps_0]$ for a fixed $\eps_0$.
The term \ef{99} produces a new discrepancy in the problem \ef{53} which is defined by
\bea\label{X2}
\big( \nabla \mathfrak{u}_{j}^{\eps\pm} , \psi_{j}^{\pm}\nabla \mathfrak{X}_j^\pm 
\big)_{\Pi^\eps} - \big( \mathfrak{u}_{j}^{\eps\pm} \nabla \mathfrak{X}_j^\pm , \nabla \psi_{j}^{\pm}
\big)_{\Pi^\eps} -\sfk^2 \big(  \mathfrak{X}_j^\pm\mathfrak{u}_{j}^{\eps\pm} , \psi_{j}^{\pm}
\big)_{\Pi^\eps}.
\eea
To estimate \ef{X2} properly, we fix $\tau > \gamma$ and $\sigma$
satisfying \ef{98}. The support of $\nabla \mathfrak{X}_j^\pm$ is contained in the union of the sets
\beas
\Sigma_{j \pm}^0 := \overline{\Pi} \cap \big( \overline{\bbB_{\widetilde{R}^{\pm}_j }(P_j^\pm)}
\setminus \bbB_{R^{\pm}_j}(P_j^\pm) \big) ,\qquad\qquad
\Sigma_{j \pm}^\eps := \overline{\varpi_j^\eps} \cap\big( \overline{\bbB_{\widetilde{R}^{\pm}_j }(P_j^\pm)}
\setminus \bbB_{R^{\pm}_j}(P_j^\pm) \big).
\eeas
In $\Sigma_{j \pm}^0$, the weights $\rho_{-\gamma , p}(x)$, $p=0,1$, are of order $1$ whereas in $\Sigma_{j \pm}^\eps$, they are of order $\eps^{-\gamma}$. In view of \ef{FeN}, the moduli of the first two terms in \ef{X2} does not exceed
\begin{equation}\label{X3}
\begin{array}{ll}
& c \,\big( \Vert \mathfrak{u}_{j}^{\eps\pm} ; H^1(\Sigma_{j \pm}^0) \Vert
+ \eps^{\gamma} \Vert \mathfrak{u}_{j}^{\eps\pm} ; H^1(\Sigma_{j \pm}^\eps) \Vert \big) \Vert \psi; W_{-\beta, -\gamma}^1(\Pi^\eps) \Vert\\[6pt]
\le & c \,\big( \eps^\tau  \Vert \cU_{j}^{\eps\pm} ; \cW_{\tau, \sigma}^1
(\bbR_+^2) \Vert  + \eps^{\gamma} e^{- \sigma R^{\pm}_j / \eps} 
 \Vert \cU_{j}^{\eps\pm} ; \cW_{\tau, \sigma}^1 (\theta_{j}^\pm\setminus\bbR_+^2 ) \Vert 
\big)\Vert \psi; W_{-\beta, -\gamma}^1(\Pi^\eps) \Vert \\[6pt]
 \le &c \eps^{\tau - \gamma} N^\eps
\Vert \psi; W_{-\beta, -\gamma}^1(\Pi^\eps) \Vert .
\end{array}
\end{equation}
The last term in \ef{X2} is estimated as follows:
\begin{equation}\label{X4}
\begin{array}{ll}
& \big|( \mathfrak{X}_j^\pm \mathfrak{u}_{j}^{\eps\pm} ,\psi_{j}^{\pm})_{
\Pi^\eps} \big| \\[6pt]
\le & c \big \Vert  (\rho_{-\gamma, 0}^\eps)^{-1} 
\mathfrak{u}_{j}^{\eps\pm} ; L^2(\Pi^\eps \cap \bbB_{\widetilde{R}^{\pm}_j} ( P_j^\pm)) 
\Vert \, \Vert \psi_{j}^{\pm} ;
W_{-\beta, -\gamma}^1(\Pi^\eps) \Vert\\[6pt]
\le &
c \eps^{1 + \gamma } \big( 
\Vert  ( 1 + |\xi|)^{\gamma -1} \cU_{j}^{\eps\pm} ; L^2(\bbR_+^2) \Vert+\Vert  ( 1 + |\eta_j^\pm|)\cU_{j}^{\eps\pm} ; L^2(\theta_{j}^\pm\setminus\bbR_+^2) \Vert
\big) \Vert \psi_{j}^{\pm} ; W_{-\beta, -\gamma}^1(\Pi^\eps)\Vert\\[6pt]
\le & c \eps N^\eps \Vert \psi_{j}^{\pm}; W_{-\beta, - \gamma}^1 (\Pi^\eps) \Vert.  
\end{array}
\end{equation}
Finally, gathering \ef{94}, \ef{X3} and \ef{X4}, we deduce that for all $F^{\eps}\in W_{-\beta,- \gamma}^{1,\eps}  (\Pi^\eps)^*$, we have the estimate
\[
\|({\rm Id} - \bfA_{\beta, \gamma}^{\eps}( \lambda) 
\bfR_{\beta, \gamma}^{\eps}  ( \lambda)) F^{\eps}; W_{-\beta,- \gamma}^{1,\eps}  (\Pi^\eps)^* \|\le c\,(\eps^{\gamma+1/2} + \eps^{\tau - \gamma} + \eps)\,\| F^{\eps}; W_{-\beta,- \gamma}^{1,\eps}  (\Pi^\eps)^*\|.
\]
This completes the proof of Theorem \ref{th4E} showing in particular that when $\gamma\in(-1/2;0)$, the operator $\bfA_{\beta, \gamma}^{\eps}( \lambda)$ is invertible for $\eps$ small enough.

\subsection{Derivation of the error estimate.} \label{sec4.5}
Now we prove that the function $u_{\rm as}^\eps$ defined in \ef{G1}--\ef{G2} yields a good approximation of $u^{\eps}_-$ as $\eps$ goes to zero. This will give us directly the proof of Theorem \ref{FF}. Estimating the discrepancies left by the asymptotic solution 
$u_{\rm as}^\eps$ is much simpler than in \S\ref{sec4.4} because the terms of the representations \ef{G1}--\ef{G2} are smooth and because we can deal with a functional $f(v) = (f,v)_{\Pi^\eps}$ which is continuous in the $L^2$-norm. The boundary conditions \ef{5} on $\Gamma_0$ and \ef{4} on 
$\Gamma_d \setminus \cup_\pm \cup_ {j=1}^J \overline{\Xi_j^{
\eps \pm}}$ are satisfied due to our choice of cut-off functions. 
Furthermore, in $\Pi$ we have
\begin{equation}\label{X5}
\begin{array}{lcl}
 - \Delta u_{\rm as}^\eps +\sfk^2 u_{\rm as}^\eps &=&\dsp
- [\Delta, X_0^\eps] \Big( w^+ - \sum_\pm \sum_{j=1}^J 
\cX_j^\pm \big( w^+ (P_j^\pm) - c_j^{2 \pm} (y - y^{\pm}_{j} )\big) \Big) \\
&&\dsp-\eps [\Delta, X_0^\eps] \Big( u_-'- \sum_\pm \sum_{j=1}^j
\cX_j^\pm Q_j^\pm \Big)\\
&&\dsp- \eps\sum_\pm \sum_{j=1}^J [\Delta, \cX_j^\pm] 
\big( \cU_j'^\pm -Q_j^\pm - c_j^{2 \pm} \xi_{j1}^\pm \big)\\ 
&&\dsp+\sfk^2 \sum_\pm \sum_{j=1}^J \cX_j^\pm \big( (1- X_0^\eps) \cU_j^{0\pm}
+ \eps \big( \cU_j'^\pm - X_0^\eps (Q_j^\pm + c_j^{2 \pm} \xi_{j1}^\pm ) \big)  \\
&=:&\dsp I_0^\eps + \eps {I_0^\eps} '
+ \sum_\pm \sum_{j=1}^J (\eps {I_{j\pm}^\eps} \!\! ' 
+\sfk^2 I_{j \pm}^\eps ), 
\end{array}
\end{equation}
where  $Q_j^\pm(x) = c_j^{0 \pm} (\ln \eps) - c_j^{1 \pm} \pi^{-1}
\ln |\eps^{-1} (x-P_j^\pm)|$. Above, the commutator $[\Delta, X_0^\eps](\cdot)$ is defined by $[\Delta, X_0^\eps](\cdot)=\Delta(X_0^\eps\cdot)-X_0^\eps\Delta\cdot$. All terms in \ef{X5} have compact supports, and in order
to estimate the scalar product 
\beas
( -\Delta u_{\rm as}^\eps +\sfk^2 u_{\rm as}^\eps, 
\psi
)_\Pi \ \ \mbox{with} \ \psi \in W_{-\beta, -\gamma}^{1, \eps}
(\Pi^\eps),
\eeas
we need to evaluate the norms $\Vert (\eps + \rho_0 )^{1 +\gamma}
I_{\ldots}^{\ldots} ; L^2(\Pi)\Vert$, see \ef{75}, \ef{76}. Since
\beas
I_0^\eps = - [\Delta, X_0^\eps ] \widetilde{\widetilde{w}}^+
, \ \ \mbox{where} \
\widetilde{\widetilde{w}}^+(x)  = O(|x - P_j^\pm|^2) \ \mbox{as} \ x \to P_j^\pm
\eeas
and supp\,$I_0^\eps = \overline{\bbR_+^2}  \cap 
\big(\cup_{\pm}\cup_{j=1}^J\overline{\bbB_{ \eps\widetilde{R}^{\pm}_j} ( P_j^\pm})\big)$, we have
\beas
\Vert (\eps + \rho_0)^{1 + \gamma} I_0^\eps ; L^2(\Pi)\Vert
\leq c \eps^{\gamma +2} .
\eeas
The same bound holds for $\Vert ( \eps + \rho_0)^{1 + \gamma}
\eps {I_0^\eps} ' ; L^2(\Pi)\Vert$ because we have
\beas
u_-'(x) - Q_j^\pm(\eps^{-1} (x-P_j^\pm) ) = 
O(|x - P_j^\pm|) \qquad \mbox{as} \ x \to P_j^\pm.
\eeas
Concerning the term $\eps {I_{j\pm}^\eps} \!\! ' $, we observe
that supp\,${I_{j\pm}^\eps} \!\! ' \subset \overline{\Pi}  \cap 
\big( \overline{\bbB_{\widetilde{R}_j^{\pm}} ( P_j^\pm)} \setminus
\bbB_{ R_j^{\pm}} ( P_j^\pm) \big) $
and 
\beas
\cU_j'^\pm (\xi_j^\pm) - Q_j^\pm(\xi_j^\pm) -c_j^{2 \pm} \xi_{j1}^\pm 
= O(|\xi|^{-1} ) \qquad\mbox{as}\ |\xi| \to +\infty, \ \xi \in \bbR_+^2 .
\eeas
Hence, we have $\Vert (\eps + \rho_0)^{1 + \gamma} \eps 
{I_{j\pm}^\eps} \!\! '  ; L^2(\Pi) \Vert \leq c\eps^2$. Finally, we obtain
\[
 \Vert (\eps + \rho_0)^{1 + \gamma}{I_{j\pm}^\eps}  
; L^2(\Pi) \Vert \le c\,\Big( \int_0^{\eps \widetilde{R}_j^{\pm}} 
( \eps +r)^{2(1 + \gamma)} r dr + 
\eps^2\int_0^{\widetilde{R}_j^{\pm}}  ( \eps +r)^{2(1+\gamma)} 
\frac{r dr}{(1 + r/ \eps)^2} \Big)^{1/2} \leq c \eps^{\gamma+2}.
\]
Let us consider the function \ef{G2} in the tunnel $\varpi_{j}^\eps$.
First of all we notice that the boundary conditions \ef{4} 
are satisfied on $\partial \Pi^\eps \cap {\rm supp} \, \cX_j^\pm$
because of our assumption on the straight segments of $\partial
\varpi_{j}^\pm$. Moreover, the harmonic functions
$\cU_j^\pm = w^\pm (P_j^\pm)$ and $\cU_j'^\pm$ are constants, see \ef{E2},
and differ from a linear function by an exponentially decaying remainder in 
$\theta_{j}^\pm\setminus\bbR_+^2$. We have
\begin{equation}\label{X6}
\begin{array}{lcl}
- \Delta u_{\rm as}^\eps +\sfk^2 u_{\rm as}^\eps&=&\dsp- X_j^\eps (\Delta -k^2)  \Big(v_j^0 + \eps 
\sum_\pm \cX_j^\pm c_{j \pm}'(\ln \eps) \Big)\\
&&\dsp-[\Delta, X_j^\eps] \Big( v_j^0 - 
\sum_\pm \cX^{\pm}_j \big( v_j^0 (\pm \ell_j) + (s \mp \ell_j) \partial_s
v_j^0(\pm \ell_j) \big) \Big)\\
&&\dsp-\eps\sum_\pm [\Delta, \cX_j^\pm] 
\big( \cU_j'^\pm - c_{j  \pm}'(\ln \eps) \mp
\eta_j^\pm \partial_s v_j^0 (\pm \ell_j)  \big) \\
&&\dsp+\sfk^2 \sum_\pm  \cX_j^\pm \big( (1- X_j^\eps) \cU_j^\pm
+ \eps \big( \cU_j'^\pm - c_{j\pm}'(\ln \eps)
\mp \eta_j^\pm  \partial_s v_j^0 (\pm\ell_j)  \big) 
\\
&=:& I_j^\eps + I_{j0}^\eps+ \eps I_j'^\eps
+\sfk^2 I_{j \pm}^\eps . 
\end{array}
\end{equation}
Because of the weight $\rho_{-\gamma,0}^\eps =
(\eps)^{-\gamma} (\eps + \rho_j)^{-1}$ in the norm
of the function $\psi \in W_{-\beta, -\gamma}^{1 , \eps} 
(\Pi^\eps)$, we need to process the norms 
$\eps^\gamma \Vert ( \eps + \rho_j) I_{\ldots}^{\ldots} \, ; L^2
(\varpi_{j}^\eps) \Vert $ of the expressions in \ef{X6}. Owing to the Taylor formula \ef{44} we have
\beas
\eps^\gamma \Vert ( \eps + \rho_j) I_{j0}^\eps ; 
L^2(\varpi_{j}^\eps) \Vert \leq c \eps^{2 + \gamma } .
\eeas
Recalling the above-mentioned exponential decay yields
\beas
\eps^{\gamma +1} \Vert ( \eps + \rho_j) 
I_{j}'^\eps ;  L^2(\varpi_{j}^\eps) \Vert 
\leq c e^{-\delta/ \eps}, \ \ \delta > 0. 
\eeas
Moreover, 
\beas
\eps^\gamma \Vert ( \eps + \rho_j) I_{j\pm}^\eps ; 
L^2(\varpi_{j}^\eps) \Vert \le c \eps^{ \gamma + 1/2  } 
\Big( \int_0^{\eps \widetilde{R}_j^{\pm}} (\eps+ \rho)^2 d \rho + 
\eps^2\int_0^{\widetilde{R}_j^{\pm} } (\eps+ \rho)^2 e^{-2\delta \rho / \eps}
 d \rho \Big)^{1/2} \leq c \eps^{2 +\gamma } .
\eeas
Here, we have taken into account the width of the tunnel $O(\eps)$ and
the evident relation 
\beas
\int_0^{\widetilde{R}_j^{\pm}} \rho^t e^{-2 \delta \rho / \eps} d \rho
\leq c_{t, \delta}\,\eps^{t+1}. 
\eeas
For the remaining term $I_j^\eps$, proceeding as in Section \ref{sec4.4}, Step 3, we obtain 
\begin{equation}\label{Question}
\eps^\gamma \Vert ( \eps + \rho_j) I_j^\eps ; 
L^2(\varpi_{j}^\eps) \Vert \leq c \eps^{\gamma+3/2}(1+|\ln\eps|).
\end{equation}
Comparing with the other previous estimates, we find that this is the largest term. This concludes the proof of Theorem \ref{FF}. Indeed, the positive number $\delta = - \gamma$ in \ef{F5}, cf. \ef{64}, can be chosen as small as one wishes. However, it is not possible to set $\delta = 0$ because $R_-''(\ln \eps)$ is a linear function of $\ln \eps$.

\subsection{Fine-tuning.}\label{sec5.3}
Now we explain how to prove that there is $\varrho > 0$ small enough such that the constant $c_{\delta}$ appearing in the estimate (\ref{F5}) of Theorem \ref{FF} can be chosen independent of $h\in\overline{\bbB_\varrho}= \{ h \in \bbR^2 : |h| \leq \varrho \}$ for all $\eps\in(0;\eps_0(\varrho)]$. Let 
\bea
\Pi^\eps(h) = \Pi \cup  \varpi_1^{\eps}(h) \cup 
\varpi_2^{\eps}(h) \cup \varpi_3^\eps \cup
\ldots \cup \varpi_J^\eps \label{R5}
\eea
be the channel \ef{2} with two tunnels shifted according to \ef{R3}. We set $x^h = (y^h,z)$ with 
\bea\label{R6}
y^h = \big( 1 - \chi_1(y) -\chi_2(y) \big) y + \chi_1(y)(y-h_1) +\chi_2(y)
(y- h_2)
\eea
and $\chi_j$ is a smooth cut-off function such that
\beas
& & \chi_j(y) = 1 \ \mbox{for} \ 
y \in \Big( y_j^- -\frac13(y_j^- - y_{j-1}^+) , y_j^+ + \frac13(
y_{j+1}^- - y_{j}^+ \Big) , \row
\chi_j(y) = 0 \ \mbox{for} \ 
y \notin \Big( y_j^- -\frac23(y_j^- - y_{j-1}^+ ), y_j^+ + \frac23(
y_{j+1}^- - y_{j}^+ \Big) . 
\eeas
The change of coordinates $x \mapsto x^h$ is non-singular for small $h \in \bbR^2$ and transforms $\Pi^\eps(h)$ into $\Pi^\eps$. Moreover, it is independent of $\eps$. The corresponding 
change affects the  differential operators of the problem \ef{3}--\ef{5} 
only in the Laplacian $\Delta$ which turns into 
\beas
\Delta + \mathscr{D}(h; x ,\nabla_x) ,
\eeas
where $\mathscr{D}$ is a first order differential operator with smooth and compactly supported coefficients, depending smoothly also on the parameter $h=(h_1 ,h_2)$. Then working as in the classical proofs of perturbations theory for linear operators (see \cite[Chapter 7, \S6.5]{Kato95}, \cite[Chap. 4]{HiPh57}), we can use the above change of variables to compare the solutions of the problem \ef{3}--\ef{5} in the same geometry $\Pi^{\eps}$ and show that they have smooth dependence with respect to $h$. We refer the reader to \cite[\S6.3]{na648} for more details.

\subsection{Well-posedness for the corrector problem}
Here we give the proof of a result of well-posedness used in the construction of the asymptotic expansions. 
\begin{proposition}\label{PropositionExistenceTerm1}
The problem \ef{My27} admits solutions $\cY^{1\pm}$, $\cY^{2\pm}$ in $H^1_{\mrm{loc}}(\overline{\Theta^{\pm}})$ with the expansions \ef{30} and \ef{30N}.   
\end{proposition}
\noindent Proof. First, for a source term $\cF$ satisfying the conditions \ef{ConstraintSource} below, we consider the problem
\begin{equation}\label{27}
\begin{array}{|rlcl}
\  -\Delta  \cU^{\pm} &=& \cF& \quad\mbox{ in } \Theta^{\pm} \\[3pt]
\ \partial_{\nu} \cU^{\pm} &=& 0&\quad\mbox{ on } \partial\Theta^{\pm} 
\end{array}
\end{equation}
The natural variational formulation of the Neumann problem \ef{27}, cf. \cite{Lad}, writes
\bea\label{27N}
\big( \nabla \cU ,\nabla \cV\big)_{\Theta^{\pm}}= ( \cF, \cV)_{\Theta^{\pm}} , \qquad \forall \cV \in \cH. 
\eea
Here $( \cdot , \cdot )_{\Theta^{\pm}}$ is the scalar  product of $L^2(\Theta^{\pm})$
and $\cH$ is the completion of $\mathscr{C}_c^\infty( \overline{ \Theta^{\pm}})$ (infinitely
differentiable and compactly supported functions) with respect to the norm
\beas
\Vert \cV; \cH \Vert := 
\big( \Vert \nabla \cV ; L^2 (\Theta^{\pm}) \Vert^2 +
\Vert \cV ; L^2 (\Theta^{\pm} \cap \bbB_R(0) ) \Vert^2 \big)^{1/2} .
\eeas
Note that the constant functions are contained in $\cH$. Using one-dimensional Hardy inequalities, one can prove that this norm is equivalent with the weighted norm
\beas
 \big( \Vert \nabla \cV ; L^2(\Theta^{\pm}) \Vert^2 
+ \Vert (1+ |\xi|)^{-1}  \cV ; L^2(\theta^{\pm}_\top) \Vert^2 +
\Vert (1+ |\xi|)^{-1} (2 + \ln|\xi|)^{-1}  \cV ; L^2(\bbR_+^2) \Vert^2 
\big)^{1/2} ,
\eeas
where $\theta^{\pm}_\top := \theta^{\pm} \setminus \overline{\bbR^2_+}$.  As a consequence, the only solution of the homogeneous problem \ef{27N} (with $\cF= 0$ ) is the constant solution in $\cH$. Moreover, according to the Fredholm theory, given a right-hand side $\cF \in L_{\rm loc}^2 (\Theta^{\pm})$ such that 
\begin{equation}\label{ConstraintSource}
(1+ |\xi|) \cF \in L^2(\theta^{\pm}_\top)\qquad\mbox{ and } \qquad
(1+ |\xi|) (2 + \ln |\xi|) \cF \in L^2(\bbR_+^2),
\end{equation}
a solution to \ef{27N} exists provided that $\cF$ satisfies the compatibility condition 
\bea\label{F10}
(\cF ,1)_{\Theta^{\pm}} = 0.
\eea
Now we show the existence of $\cY^{1\pm}$, $\cY^{2\pm}$ in $H^1_{\mrm{loc}}(\overline{\Theta^{\pm}})$ solving \ef{My27} with the expansions \ef{30} and \ef{30N}. Introduce a cut-off function $\chi\in\mathscr{C}^{\infty}(\overline{\Theta^{\pm}})$ such that $\chi=0$  in $\bbB_R(0)$ and  $\chi=1$  in $\bbB_{2R}(0)$ for $R$ large enough. Define the functions $\aleph^{1\pm}$, $\aleph^{2\pm}$ such that 
\[
 \aleph^{1\pm} (\xi) = \left\{
\begin{array}{ll}
 \dsp\frac{\chi}{\pi} \ln \frac{1}{|\xi|} , &\qquad \mbox{in }\bbR_+^2 \\[7pt]
\chi\dsp\frac{\eta^{\pm}}{H(\pm\ell)} &\qquad\mbox{in }\theta^{\pm}\setminus\bbR_+^2 
\end{array}
\right.\qquad\mbox{ and }\qquad\aleph^{2\pm} (\xi) = \left\{
\begin{array}{ll}
\dsp\chi\,\xi_1 & \qquad\mbox{in }\bbR_+^2 \\[7pt]
\dsp 0 &\qquad\mbox{in }\theta^{\pm}\setminus\bbR_+^2. 
\end{array}
\right.
\]
Set $\cF^{1\pm}=\Delta\aleph^{1\pm}$ and $\cF^{2\pm}=\Delta\aleph^{2\pm}$. Observe that these functions are compactly supported. Moreover, a direct calculus shows that $\cF^{1\pm}$, $\cF^{2\pm}$ satisfy the compatibility condition \ef{F10}. Denote $\cV^{1\pm}$, $\cV^{2\pm}$ the corresponding solutions of \ef{27} which behave as $O(|\xi|^{-1})$ at infinity in $\bbR_+^2$ (remember that the solution is defined up to an additional constant). Finally set $\cY^{1\pm}:=\cV^{1\pm}+\aleph^{1\pm}$ and $\cY^{2\pm}:=\cV^{2\pm}+\aleph^{2\pm}$. Using Fourier decomposition one can verify that $\cY^{1\pm}$, $\cY^{2\pm}$ admit the desired behaviours at infinity. \ \ $\boxtimes$

\section{Conclusion}

In this article, we presented a method to construct non-reflecting perturbations of the bottom of a channel for a water-wave problem. To proceed, we considered singular perturbations of the geometry with thin curved channels. With this approach, we showed how to get $R=0$ and we proved that we cannot get $T=1$. As a consequence, the transmitted field exhibits a phase shift with respect to the incident field. Let us mention that acoustic waveguides and water channels are somehow opposite of each other in the following sense: for the former, regular perturbations of the boundary, Fig.\,\ref{fig5},\,a), may at most be non-reflecting ($R^{\eps}=0$) while wells,  Fig.\,\ref{fig5},\,b), may be perfectly invisible ($T^{\eps}=1$). For the latter these properties are reversed. To obtain the main result of this article (construction of non reflecting perturbations of the bottom), we extend the results known for asymptotic expansion in junctions of massive bodies and thin ligaments. More precisely, we considered situations where the ligaments have a non constant width and where the ligaments do not arrive perpendicularly to the massive body. These studies were not performed in literature. Finally, let us mention that another approach to achieve invisibility has been proposed in \cite{ChNP18,ChPaSu} for acoustics problems. It would be relevant to study if we can adapt it to deal with water-wave problems.

\section*{Acknowledgments} 
The second named 
author was partially supported by RFFI, grant 18-01-00325
and by the Academy of Finland grant no. 139545. The third named
author was partially supported by a research grant from Faculty of Science of the University of Helsinki.

\bibliography{Bibli}

\def\cprime{$'$}
\begin{thebibliography}{10}

\bibitem{Arse76}
{A.A.} Arsen'ev.
\newblock The existence of resonance poles and scattering resonances in the
  case of boundary conditions of the second and third kind.
\newblock {\em USSR Comput. Math. Math. Phys.}, 16(3):171--177, 1976.

\bibitem{AsPV00}
A.~Aslanyan, L.~Parnovski, and D.~Vassiliev.
\newblock Complex resonances in acoustic waveguides.
\newblock {\em Quart. J. Mech. Appl. Math.}, 53(3):429--447, 2000.

\bibitem{A4}
{F.L.} Bakharev and {S.A.} Nazarov.
\newblock Gaps in the spectrum of a waveguide composed of domains with
  different limiting dimensions.
\newblock {\em Sib. Math. J.}, 56(4):575--592, 2015.

\bibitem{Beal73}
{J.T.} Beale.
\newblock Scattering frequencies of resonators.
\newblock {\em Comm. Pure Appl. Math.}, 26(4):549--563, 1973.

\bibitem{na648}
{A.-S.} {Bonnet-Ben Dhia}, {L.} Chesnel, and {S.A.} Nazarov.
\newblock Perfect transmission invisibility for waveguides with sound hard
  walls.
\newblock {\em J. Math. Pures Appl.}, 111:79--105, 2018.

\bibitem{BN}
{A.-S.} {Bonnet-Ben Dhia} and {S.A.} Nazarov.
\newblock Obstacles in acoustic waveguides becoming ``invisible'' at given
  frequencies.
\newblock {\em Acoust. J.}, 59(6):685--692, 2013.
\newblock English transl. Acoust. Phys. 59(6):633--639, 2012.

\bibitem{na582}
{A.-S.} {Bonnet-Ben Dhia}, {S.A.} Nazarov, and {J.} Taskinen.
\newblock Underwater topography ``invisible'' for surface waves at given
  frequencies.
\newblock {\em Wave Motion}, 57(0):129--142, 2015.

\bibitem{ChNP18}
L.~Chesnel, {S.A.} Nazarov, and V.~Pagneux.
\newblock Invisibility and perfect reflectivity in waveguides with finite
  length branches.
\newblock {\em SIAM J. Appl. Math.}, 78(4):2176--2199, 2018.

\bibitem{ChPaSu}
L.~Chesnel and V.~Pagneux.
\newblock Simple examples of perfectly invisible and trapped modes in
  waveguides.
\newblock {\em Quart. J. Mech. Appl. Math.}, 71(3):297--315, 2018.

\bibitem{EvLV94}
{D.V.} Evans, M.~Levitin, and D.~Vassiliev.
\newblock Existence theorems for trapped modes.
\newblock {\em J. Fluid. Mech.}, 261:21--31, 1994.

\bibitem{Gady93}
{R.R.} Gadyl'shin.
\newblock {Characteristic frequencies of bodies with thin spikes. I.
  Convergence and estimates}.
\newblock {\em Math. Notes}, 54(6):1192--1199, 1993.

\bibitem{A24}
{R.R.} Gadyl'shin.
\newblock {On the eigenvalues of a ``dumbbell with a thin handle''.}
\newblock {\em Izv. Math.}, 69(2):265--329, 2005.

\bibitem{Have29}
{T.H.} Havelock.
\newblock Forced surface waves.
\newblock {\em Phil. Mag.}, 8:569--576, 1929.

\bibitem{HiPh57}
E.~Hille and {R.S.} Phillips.
\newblock {\em Functional analysis and semi-groups}, volume~31.
\newblock Amer. Math. Soc., 1957.

\bibitem{Ilin}
A.~M. Il'in.
\newblock {\em Matching of asymptotic expansions of solutions of boundary value
  problems}, volume 102 of {\em Translations of Mathematical Monographs}.
\newblock AMS, Providence, RI, 1992.

\bibitem{A30}
P.~Joly and S.~Tordeux.
\newblock Asymptotic analysis of an approximate model for time harmonic waves
  in media with thin slots.
\newblock {\em Math. Mod. Num. Anal.}, 40(1):63--97, 2006.

\bibitem{A31}
P.~Joly and S.~Tordeux.
\newblock {Matching of asymptotic expansions for wave propagation in media with
  thin slots I: The asymptotic expansion}.
\newblock {\em SIAM Multiscale Model. Simul.}, 5(1):304--336, 2006.

\bibitem{A32}
P.~Joly and S.~Tordeux.
\newblock {Matching of asymptotic expansions for waves propagation in media
  with thin slots II: The error estimates}.
\newblock {\em Math. Mod. Num. Anal.}, 42(2):193--221, 2008.

\bibitem{Kato95}
T.~Kato.
\newblock {\em {Perturbation theory for linear operators.}}
\newblock Springer-Verlag, Berlin, reprint of the corr. print. of the 2nd ed.
  1980 edition, 1995.

\bibitem{Ko67}
V.~A. {Kondratiev}.
\newblock Boundary-value problems for elliptic equations in domains with
  conical or angular points.
\newblock {\em Trudy Moskov. Matem. Obshch.}, 16:209--292, 1967.
\newblock English transl. Trans. Moscow Math. Soc. 16:227--313, 1967.

\bibitem{KoMR01}
V.~A. Kozlov, V.~G. Maz'ya, and J.~Rossmann.
\newblock {\em Spectral problems associated with corner singularities of
  solutions to elliptic equations}, volume~85 of {\em Mathematical Surveys and
  Monographs}.
\newblock AMS, Providence, 2001.

\bibitem{KoMM94}
{V.A.} Kozlov, {V.G.} Maz'ya, and {A.B.} Movchan.
\newblock Asymptotic analysis of a mixed boundary value problem in a
  multi-structure.
\newblock {\em Asymptot. Anal.}, 8(2):105--143, 1994.

\bibitem{KMM2}
{V.A.} Kozlov, {V.G.} Maz'ya, and {A.B.} Movchan.
\newblock Fields in non-degenerate 1d-3d elastic multi-structures.
\newblock {\em Quart. J. Mech. Appl. Math.}, 54(2):177--212, 2001.

\bibitem{KuMaVa}
N.~Kuznetsov, {V.G.} Maz'ya, and B.~Vainberg.
\newblock {\em Linear water waves: a mathematical approach}.
\newblock Cambridge University Press, 2002.

\bibitem{Lad}
{O.A.} Ladyzhenskaya.
\newblock Boundary value problems of mathematical physics.
\newblock {\em Nauka, Moscow}, 1973.
\newblock English transl. Springer, New York, 1985.

\bibitem{LeKi01}
{H.-W.} Lee and {C.S.} Kim.
\newblock Effects of symmetries on single-channel systems: Perfect transmission
  and reflection.
\newblock {\em Phys. Rev. B}, 63(7):075306, 2001.

\bibitem{Lion06}
{J.-L.} Lions.
\newblock {\em Perturbations singuli{\`e}res dans les probl{\`e}mes aux limites
  et en contr{\^o}le optimal}, volume 323.
\newblock Springer, 2006 (from the original version of 1973).

\bibitem{MaNaPl}
{V.G.} {Maz'ya}, {S.A.} Nazarov, and {B.A.} Plamenevski{\u\i}.
\newblock {\em {Asymptotic theory of elliptic boundary value problems in
  singularly perturbed domains, Vol. 1}}.
\newblock {Birkh\"{a}user}, Basel, 2000.
\newblock Translated from the original German 1991 edition.

\bibitem{MrMK11}
A.E. Miroshnichenko, B.A. Malomed, and {Y.S.} Kivshar.
\newblock Nonlinearly {PT}-symmetric systems: Spontaneous symmetry breaking and
  transmission resonances.
\newblock {\em Phys. Rev. A}, 84(1):012123, 2011.

\bibitem{Naza96}
{S.A.} Nazarov.
\newblock Junctions of singularly degenerating domains with different limit
  dimensions 1.
\newblock {\em J. Math. Sci. (N.Y.)}, 80(5):1989--2034, 1996.

\bibitem{na285}
{S.A.} Nazarov.
\newblock Junctions of singularly degenerating domains with different limit
  dimensions. 2.
\newblock {\em Trudy seminar. Petrovskii., Moscow Univ.}, 20:155--195, 1997.
\newblock English transl. J. Math. Sci., 97(3):155--195, 1999.

\bibitem{na262}
{S.A.} Nazarov.
\newblock The polynomial property of self-adjoint elliptic boundary-value
  problems and an algebraic description of their attributes.
\newblock {\em Russ. Math. Surv.}, 54(5):947--1014, 1999.

\bibitem{na345}
{S.A.} Nazarov.
\newblock Elliptic boundary value problems on hybrid domains.
\newblock {\em Funkt. Anal. i Prilozhen}, 38(4):55--72, 2004.
\newblock English transl. Funct. Anal. Appl., 38(4):283--297, 2004.

\bibitem{Naza05}
{S.A.} Nazarov.
\newblock Asymptotic analysis and modeling of the jointing of a massive body
  with thin rods.
\newblock {\em J. Math. Sci. (N.Y.)}, 127(5):2192--2262, 2005.

\bibitem{na489}
{S.A.} Nazarov.
\newblock Asymptotic expansions of eigenvalues in the continuous spectrum of a
  regularly perturbed quantum waveguide.
\newblock {\em Theor. Math. Phys.}, 167(2):606--627, 2011.

\bibitem{na514}
{S.A.} Nazarov.
\newblock Asymptotics of solutions to the spectral elasticity problem for a
  spatial body with a thin coupler.
\newblock {\em Sibirsk. Mat. Zh.}, 53(2):345--364, 2012.
\newblock English transl. Sib. Math. J., 53(2):274--290, 2012.

\bibitem{na546}
{S.A.} Nazarov.
\newblock Enforced stability of a simple eigenvalue in the continuous spectrum
  of a waveguide.
\newblock {\em Funct. Anal. Appl.}, 47(3):195--209, 2013.

\bibitem{na576}
{S.A.} Nazarov.
\newblock Modeling of a singularly perturbed spectral problem by means of
  selfadjoint extensions of the operators of the limit problems.
\newblock {\em Funkt. Anal. i Prilozhen}, 49(1):31--48, 2015.
\newblock English transl. Funct. Anal. Appl., 49(1):25--39, 2015.

\bibitem{NaPl}
{S.A.} Nazarov and {B.A.} Plamenevski{\u\i}.
\newblock {\em Elliptic problems in domains with piecewise smooth boundaries},
  volume~13 of {\em Expositions in Mathematics}.
\newblock De Gruyter, Berlin, Germany, 1994.

\bibitem{na328}
{S.A.} Nazarov and J.~Soko{\l}owski.
\newblock {The topological derivative of the Dirichlet integral under formation
  of a thin bridge}.
\newblock {\em Sibirsk. Mat. Zh.}, 45(2):410--426, 2004.
\newblock English transl. Sib. Math. J. 45(2):341--355, 2004.

\bibitem{na642}
{S.A.} Nazarov and J.~Taskinen.
\newblock Radiation conditions for the linear water-wave problem in periodic
  channels.
\newblock {\em Math. Nachr.}, 290(11--12):1753--1778, 2017.

\bibitem{PoGP99}
{J.A.} Porto, {F.J.} {Garcia-Vidal}, and {J.B.} Pendry.
\newblock Transmission resonances on metallic gratings with very narrow slits.
\newblock {\em Phys. Rev. Lett.}, 83(14):2845, 1999.

\bibitem{Shao94}
{Z.-A.} Shao, W.~Porod, and {C.S.} Lent.
\newblock Transmission resonances and zeros in quantum waveguide systems with
  attached resonators.
\newblock {\em Phys. Rev. B}, 49(11):7453, 1994.

\bibitem{VD}
M.~Van~Dyke.
\newblock {\em Perturbation methods in fluid mechanics}.
\newblock The Parabolic Press, Stanford, Calif., 1975.

\bibitem{ViLu}
{L.A.} Vishik, {M.I.};~Ljusternik.
\newblock The asymptotic behaviour of solutions of linear differential
  equations with large or quickly changing coefficients and boundary
  conditions.
\newblock {\em Uspehi. Mat. Nauk.}, 15(4):27--95, 1960.
\newblock English transl. Russian Math. Surveys, 15(4):23--91, 1960.

\bibitem{Zhuk10}
{S.V.} Zhukovsky.
\newblock Perfect transmission and highly asymmetric light localization in
  photonic multilayers.
\newblock {\em Phys. Rev. A}, 81(5):053808, 2010.

\end{thebibliography}
\bibliographystyle{plain}

\end{document}